%% file: main.tex
\documentclass[journal,letterpaper]{IEEEtran}

\usepackage[numbers,sort&compress,square]{natbib}   
\usepackage{graphicx}
\usepackage{subfigure}
\usepackage{amsmath,amssymb}
\usepackage{array}
\usepackage{amsthm}
\usepackage{mathrsfs}
\usepackage{amsmath}
\usepackage{setspace}
\usepackage{color}
\usepackage{nnfootnote}
\usepackage[finalnew]{trackchanges} 

\newcommand{\beq}{\begin{equation}}
\newcommand{\eeq}{\end{equation}}
\newcommand{\bea}{\begin{eqnarray}}
\newcommand{\eea}{\end{eqnarray}}
\newcommand{\bean}{\begin{eqnarray*}}
\newcommand{\eean}{\end{eqnarray*}}

\newcommand{\expE}{\mathbb{E}}
\newcommand{\prob}{\mathbb{P}}
\newcommand{\arrival}{\mathcal{A}}
\newcommand{\calV}{\mathcal{V}}
\newcommand{\calk}{\mathcal{K}}

\newcommand{\argmin}{\operatornamewithlimits{arg\,min}}

\newcommand{\delnote}[1]{}

\newtheorem{theorem}{Theorem}
\newtheorem{proposition}[theorem]{Proposition}
\newtheorem{lemma}[theorem]{Lemma}
\newtheorem{corollary}[theorem]{Corollary}

\theoremstyle{definition}
\newtheorem{define}[theorem]{Definition}

\newtheorem{remark}[theorem]{Remark}

\begin{document}

\title{Opportunities for Network Coding:\\ To Wait or Not to Wait}

\author{\IEEEauthorblockN{Yu-Pin Hsu\IEEEauthorrefmark{1},
Navid Abedini\IEEEauthorrefmark{1},
Natarajan Gautam\IEEEauthorrefmark{2}, 
Alex Sprintson\IEEEauthorrefmark{1}, and
Srinivas Shakkottai\IEEEauthorrefmark{1}}  \\
\IEEEauthorblockA{\IEEEauthorrefmark{1}Department of Electrical \& Computer Engineering, Texas A\&M University} \\
\IEEEauthorblockA{\IEEEauthorrefmark{2}Department of Industrial \& Systems Engineering, Texas A\&M University}\\
\IEEEauthorblockA{\{yupinhsu, novid\_abed, gautam, spalex, sshakkot\}@tamu.edu}
}

\addeditor{AS}

\maketitle

\nnfoottext{An earlier version of this paper was presented at the ISIT, 20111.}

\input{abstract.tex}
\input{introduction.tex}
\input{model.tex}
\input{analysis.tex}

\input{numericals.tex}

\input{extension.tex}

\input{conclusion.tex}

\small
\bibliographystyle{IEEEtran}
\bibliography{IEEEabrv,ref}

\end{document}

%% file: abstract.tex
\begin{abstract}
\annote[AS]{It has been well established}{Rephrased the abstract.} that wireless network coding can significantly improve the efficiency of multi-hop wireless networks. However, in a stochastic environment some of the packets might not have coding pairs, which limits the number of available coding opportunities. In this context, an important decision is whether to delay packet transmission in hope that a coding pair will be available in the future or transmit a packet without coding. The paper addresses this problem by  formulating a stochastic dynamic program whose objective is to minimize the long-run average cost per unit time incurred due to transmissions and delays. In particular, we identify optimal control actions that would balance between costs of transmission against the costs incurred due to the delays. Moreover, we seek to address a crucial question: what should be observed as the state of the system? We analytically show that observing queue lengths suffices if the system can be modeled as a Markov decision process. We also show that a stationary threshold type policy based on queue lengths is optimal. We further substantiate our results with simulation experiments for more generalized settings.
\end{abstract}

%% file: introduction.tex

\section{Introduction}\label{sec:intro}

\note[AS]{The paper uses quite a heavy notation - I think that it would be quite beneficial to use a notation table that will include all the notation used in the paper}

\begin{figure}
\centering
\includegraphics[width=.32\textwidth]{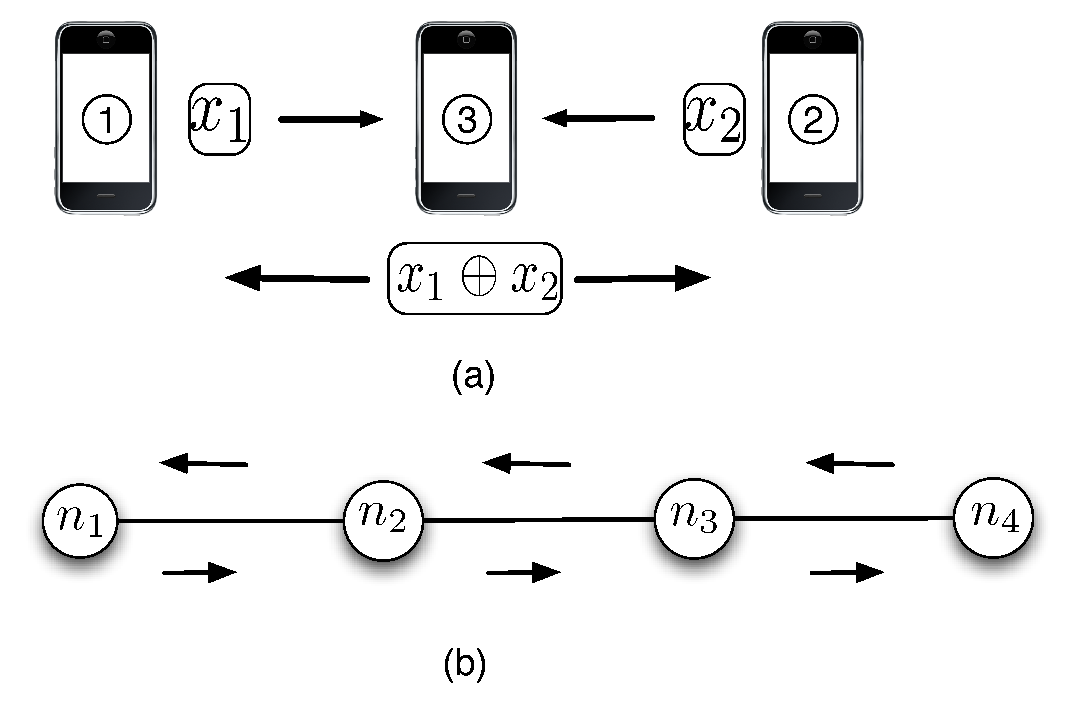}
\caption{(a) Wireless Network Coding (b) Reverse carpooling.}
\label{fig:network}
\end{figure}

In recent years, there has been a growing interest in the applications of network coding techniques in wireless  networks. It was shown that  network coding can result in significant improvements in the performance in terms of delay and transmission count. For example, consider a wireless network coding scheme depicted in Fig. \ref{fig:network}(a). Here, wireless nodes 1 and 2 need to exchange packets $x_1$ and $x_2$ through a relay node (node 3).  A simple \emph{store-and-forward} approach needs four transmissions.  In contrast, the network coding solution uses a  \emph{store-code-and-forward} approach in which the two packets $x_1$ and $x_2$ are combined by means of a bitwise XOR operation at the relay and are broadcast to nodes 1 and 2 simultaneously.  Nodes 1 and 2 can then decode this coded packet to obtain the packets they need.

Effros et al. \cite{reverse-car-pooling} introduced the strategy of reverse carpooling that allows two information flows traveling in opposite directions to share a path. Fig. \ref{fig:network}(b) shows an example of two connections, from $n_1$ to $n_4$ and from $n_4$ to $n_1$ that share a common path $(n_1,n_2,n_3,n_4)$. The wireless network coding approach results in a significant (up to 50\%) reduction in the number of transmissions for two connections that use reverse carpooling. In particular, once the first connection is established, the second connection (of the same rate) can be established in the opposite direction with little additional cost.

In this paper, we focus on the design and analysis of scheduling protocols that exploit the fundamental trade-off between the number of transmissions and delay in the reverse carpooling schemes. In particular, to cater to  delay-sensitive  applications, the network must be aware that savings achieved by coding may be offset by delays incurred in waiting for such opportunities.  Accordingly, we design delay-aware controllers that use local information to decide whether or not to wait for a coding opportunity, or to go ahead with an uncoded transmission. By sending uncoded packets we do not take advantage of network coding, resulting in a penalty in terms of transmission count, and, as a result, energy-inefficiency. However, by waiting for a coding opportunity, we might be able to achieve energy efficiency at the cost of a small delay increase.

\begin{figure}[ht]
\centering
\includegraphics[width=.35\textwidth]{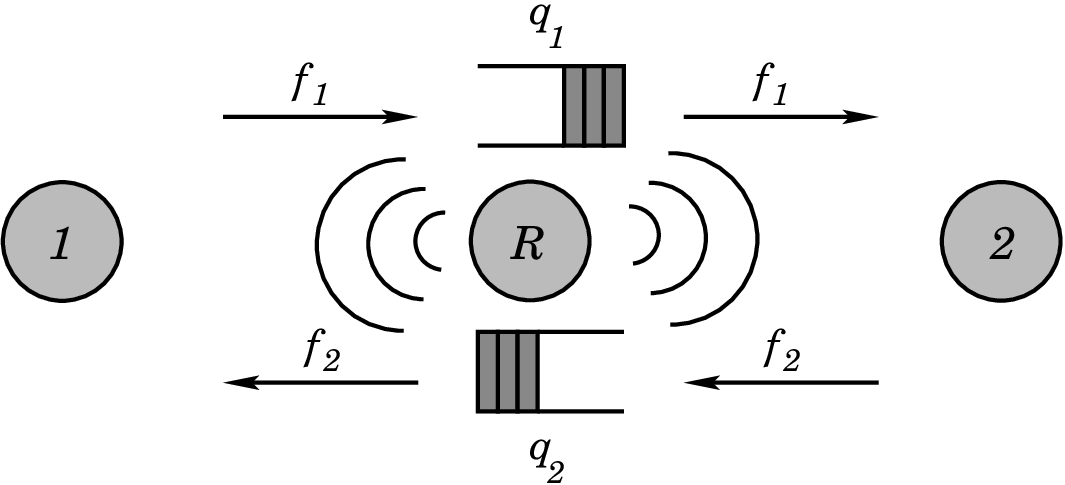}
\caption{3-Node Relay Network.}
\label{fig:relaynetwork}
\end{figure}

Consider a relay node that transmits packets between two of its adjacent nodes with flows in opposite directions, as depicted in Fig. \ref{fig:relaynetwork}. The relay maintains two queues $q_1$ and $q_2$, such that $q_1$ and $q_2$ store packets that need to be delivered to node 2 and node 1, respectively. If both queues are not empty, then it can relay two packets from both queues by performing an XOR operation. However, what should the relay do if one of the queues has packets to transmit, while the other queue is empty? Should the relay wait for a coding opportunity or just transmit a packet from a non-empty queue without coding? This is the fundamental question we seek to answer. In essence we would like to trade off efficiently transmitting the packets against high quality of service (i.e., low delays). 

\subsection{Related Work} 
Network coding research was initiated by the seminal work  of Ahlswede et al. \cite{network-information-flow} and since then attracted major interest from the research community. Network coding technique for wireless networks has been considered by Katti et al. \cite{COPE-Katabi}. They propose an architecture, referred to as COPE, which contains a special network coding layer between the IP and MAC layers. In \cite{MORE-Chachulski}, an opportunistic routing protocol is proposed, referred to as MORE, that randomly mixes packets that belong to the same flow before forwarding them to the next hop. In addition, several works, e.g., \cite{cross-layer-opt-Ephremides,cross-lay-opt-Shroff,population-game-shakkottai,remedy-packet,dynamic-alg-tracey,distributed-alg-Yeh}, investigate the scheduling and/or routing problems in the network coding enabled networks. Sagduyu and Ephremides \cite{cross-layer-opt-Ephremides} focus on the network coding in the tandem networks and formulate related cross-layer optimization problems, while Khreishah et al. \cite{cross-lay-opt-Shroff} devise a joint coding-scheduling-rate controller when the pairwise intersession network coding is allowed.  \annote[AS]{Reddy et al.}{We cannot say ``our previous work'' the set of authors was different.} \cite{population-game-shakkottai} have showed how to design coding-aware routing controllers that would maximize coding opportunities in multihop networks. References \cite{remedy-packet} and \cite{dynamic-alg-tracey} attempt to schedule the network coding between multiple-session flows. Xi and  Yeh \cite{distributed-alg-Yeh}  propose a distributed algorithm that minimizes the transmission cost of a multicast session.  

The work of Ciftcioglu et al. \cite{wait-coding} is the most relevant to our paper. It is important to note that this work was performed independently and  analyzed a related problem from a different perspective. In particular, \cite{wait-coding}  proposed a control policy that strikes a balance between the delay and the cost, as well as compares the policy to one that never waits for the coding opportunity. In contrast, in our paper, we provide a \annote[AS]{provably optimal control policy}{Check the correctness of the claim} and identify its structure. 


In this paper, we consider a stochastic arrival process and address the decision problem of whether or not a packet should wait for a coding opportunity. Our objective is therefore to study the delicate trade-off between the energy consumption and the queueing delay when network coding is an option.  We use the Markov decision process (MDP) framework to model this problem and formulate a stochastic dynamic program that determines the optimal control actions in various states. While there exists a large body of literature on the analysis of MDPs (see, e.g., \cite{MDP-Puterman, MDP-Ross, MDP-Sennott, MDP-Bertsekas}),  there is no clear methodology to find optimal policies for the problems that possess the proprieties of infinite horizon, average cost optimization, and with a countably infinite state space.  Indeed, \cite{MDP-Bertsekas} remarks that it is difficult to analyze and obtain optimal policies for such problems.  The works  in \cite{mdp-countable1,mdp-countable2,mdp-countable3,stationary-policy-Sennott} contribute to the analysis of MDPs with countably infinite state space. Moreover, reference \citep{monotone-policy-survey} that surveys the recent results on the monotonic structure of optimal policy, states that while one dimensional MDP with convex cost functions has been extensively studied, limited models for multi-dimensional spaces are dealt with due to the correlations between dimensions. In many high-dimension cases, one usually directly  investigates the properties of the cost function. As we will see later, this paper poses precisely such a problem, and showing the properties of optimal solution is one of our main contributions. 

\subsection{Main Results}


We first consider the case illustrated in Fig. \ref{fig:relaynetwork}, in which we have a single relay node with two queues that contain packets traversing in opposite directions. We assume that time is slotted, and the relay can transmit at most one packet during each time slot.  We also assume that the arrivals into each queue are independent and identically distributed.  Each transmission by the relay incurs a cost, and similarly, each time slot when a packet waits in the queue has some cost.  We  would like to minimize the average sum of the two costs.  In general, we could utilize a controller that belongs to one of the following sets \citep{MDP-Puterman}:
\begin{itemize}
	\item  $\Pi^{\text{HR}}$ - a set of randomized history dependent policies;
	\item $\Pi^{\text{MR}}$ - a set of randomized Markov policies;
	\item $\Pi^{\text{SR}}$ - a set of  randomized stationary policies;
	\item $\Pi^{\text{SD}}$ - a set of  deterministic stationary policies.
\end{itemize}

It is not hard to see (as shown in \cite{MDP-Puterman}) that
\begin{equation*}
		\Pi^{\text{SD}} \subset \Pi^{\text{SR}}  \subset \Pi^{\text{MR}} \subset \Pi^{\text{HR}}.
\end{equation*}
The complexity of the algorithms increases from left to right above:  in what regime does the solution to our problem lie?  We can think of the system state as the two queue lengths. 
We find that the optimal policy is a simple queue-length threshold policy with one threshold for each queue at the relay, and whose action is simply: if a coding opportunity exists, code and transmit; else transmit a packet if the threshold for that queue is reached.  We then show how to find the optimal thresholds.  Thus, our result implies that although waiting time information might be available, we do not need to actually use it.  

We examine two general models afterward. In the first model, the service capacity of the relay is not restricted to one packet per time slot. Then, if the relay can serve a batch of packets, we find that the optimal controller is of the threshold type for one queue, when the queue length of the other queue is fixed.  Secondly, we study an arrival process with memory (Markov modulated).  Here, we discover that the optimal policy has multiple thresholds.  

We then perform a numerically study of a number of policies that are based on waiting time and queue length, waiting time only, as well as the optimal deterministic queue-length threshold policy to  indicate the potential of our approach.  We also evaluate the performance of a deterministic queue length based policy in the line network topology via simulations. 

\textbf{Contributions.} Our contributions can be summarized as follows. We introduce the problem of delay versus coding efficiency trade-off, as well as formulate it as an MDP problem and obtain the structure of the optimal policy.  It turns out that the optimal policy does not use the waiting time information. Moreover, we prove that the optimal policy is stationary and of threshold type in terms of the queue lengths, and therefore is easy to implement.  While it is easy to analyze MDPs that have a finite number of states, or involve a discounted total cost optimization with a single communicating class, our problem does not possess any of these properties. Hence, although our policy is simple, the proof is extremely intricate. Furthermore, our policy and proof techniques can be extended to other scenarios such as batched service and Markov-modulated arrival process.      

%% file: model.tex
%

\section{System Overview}\label{sec:overview}
\subsection{System model}
Consider a multi-hop wireless network operating a time-division multiplexing scheme to store and forward packets from various sources to destinations. Time is divided into \emph{slots} that are further divided into \emph{mini-slots}.  In each slot, each node is allowed to transmit in its assigned mini-slot.  Such a deterministic schedule without interference is easy to construct, e.g., see \cite{KulVis04,ShaLiu10} for a method to do so in a unit square with randomly dropped nodes, in which every node gets a transmission opportunity with finite periodicity.  



Our first focus is on the case of a single relay node of interest, which has the potential for network coding packets from flows in opposing directions.  Consider Fig. \ref{fig:relaynetwork} again.  We call the two adjacent nodes to the relay $R$ as nodes $1$ and $2$. We assume that there is a flow $f_1$ that goes from node $1$ to $2$ and another flow $f_2$ from node $2$ to $1$, both of which are through the relay under consideration. The packets from both flows are stored at separate queues, $q_1$ and $q_2$, at node $R$. Each slot is divided into several mini-slots, such that the last mini-slot is used by the relay and all other mini-slots are used by nodes $1$ and $2$. Note that the time period between transmission opportunities for the relay is precisely one slot.

The number of arrivals between consecutive slots to both flows is assumed to be independent of each other and also independent and identically distributed (i.i.d.) over time, with the random variables \annote[AS]{$\arrival_i$}{This notation is very similar to $a_t$ - the reader might be a bit confused.} for $i=1,2$ respectively. In each slot, $n$ packets arrive  at $q_i$ with the probability $\prob(\arrival_i=n)=p^{(i)}_n$ for $n \in \mathbb{N} \cup \{0\}$. Afterward, the relay gets an opportunity to transmit. Initially we assume that the relay can transmit a maximum of one packet in each time slot. 



\subsection{Markov Decision Process Model}
We use a Markov decision process (MDP) model to develop a strategy for the relay to decide its best course of action at every transmission opportunity. For $i=1,2$ and $t = 0, 1, 2, \cdots$, let $Q^{(i)}_t$ be the number of packets in $q_i$ at the $t^{\text{th}}$ time slot just before an opportunity to transmit. Let $a_t$ be the action chosen at the end of the $t^{\text{th}}$ time slot with $a_t = 0$ implying the action is to do nothing and $a_t= 1$ implying the action is to transmit. Clearly, if $Q^{(1)}_t + Q^{(2)}_t = 0$, then $a_t = 0$ because that is the only feasible action. Also, if $Q^{(1)}_t Q^{(2)}_t > 0$, then $a_t=1$ because the best option is to transmit as a coded XOR packet as it  reduces both the number of transmissions as well as latency. However, when exactly one of $Q^{(1)}_t$ and $Q^{(2)}_t$ is non-zero, it is unclear what the best action is. 

To develop a strategy for that, we first define the costs for latency and transmission. Let $C_T$ be the cost for transmitting a packet and $C_H$ be the cost of holding a packet for a length of time equal to one slot. Without loss of generality, we assume that \annote[AS]{if a packet is transmitted in the same slot that it arrived,}{Does this rely on the assumption that the packets arrive in the beginning of the slot?} its latency is zero. Also, the cost of transmitting a coded packet is the same as that of a non-coded packet. That said, our objective is to derive an optimal policy that minimizes the long-run average cost per slot. Therefore, we define the MDP$\{(Q_t,a_t), t \ge 0\}$ where $Q_t = (Q^{(1)}_t,Q^{(2)}_t)$ is the state of the system and $a_t$ the control action chosen by the relay at the $t^{\text{th}}$ slot. The state space (i.e., all possible values of $Q_t$) is the set $\{(i,j): i=0,1, \cdots; j=0, 1, \cdots\}$.


Let $C(Q_t,a_t)$ be the \textit{immediate cost} if action $a_t$ is taken at time $t$ when the system is in state $Q_t = (Q^{(1)}_t,Q^{(2)}_t)$. Therefore,
\bea
C(Q_t, a_t) = C_H ([Q^{(1)}_t-a_t]^+ + [Q^{(2)}_t-a_t]^+) + C_T a_t,
\label{eq:immediate cost}
\eea
where $[x]^+ = \max(x,0)$. The long-run average cost for some policy $\theta \in \Pi^{\text{HR}}$ is given by
\bea 
V(\theta) = \lim_{K \rightarrow \infty} \frac{1}{K+1} \expE_{\theta} \left[ \sum_{t=0}^K C(Q_t,a_t) | Q_0 = (0,0) \right], \label{eq:long-time-average-cost}
\eea
where $\expE_{\theta}$ is the expectation operator taken for the system under policy $\theta$. Notice that our initial state is an empty system, although the average cost would not depend on it.  Our goal is to characterize and obtain the \textit{average-optimal policy}, i.e., the policy that minimizes $V(\theta)$.  We first describe the probability law for our MDP and then in subsequent sections develop a methodology to obtain the average-optimal policy. 

For the MDP$\{(Q_t,a_t), t \ge 0\}$,  let $P_{a_t}(Q_t, Q_{t+1})$ be the transition probability from state $Q_t$ to $Q_{t+1}$ associated with action $a_t\in\{0,1\}$. Then the probability law can be derived as $P_{0}\left((i,j), (k,l) \right) = p^{(1)}_{k-i} p^{(2)}_{l-j}$ for all $k \ge i$ and $l \ge j$; otherwise, $P_{0}\left((i,j), (k,l) \right) = 0$.  Also, $P_{1}\left((i,j), (k,l) \right) = p^{(1)}_{k-[i-1]^+} p^{(2)}_{l-[j-1]^+}$ for all $k \ge [i-1]^+$ and $l \ge [j-1]^+$; otherwise, $P_{1}\left((i,j), (k,l) \right) = 0$.

A list of important notation used in this paper is summarized in Table \ref{table:notation}.

\begin{table}
\begin{tabular}{|c|p{6cm}|}
\hline
$\arrival_i$ & Random variable that represents the number of packets that arrives at $q_i$ for each time slot\\
\hline
$p^{(i)}_n$ & Probability that $n$ packets arrive at  $q_i$, i.e., $\mathbb{P}(\arrival_i=n)$\\
\hline
$Q^{(i)}_t$ & The number of packets in $q_i$ at time $t$\\
\hline
$Q_t$ & System state, i.e., $(Q^{(1)}_t,Q^{(2)}_t)$ \\
\hline
$a_t$ & Action chosen by relay at time $t$\\
\hline
$C_T$ & Cost of transmitting one packet\\
\hline
$C_H$ & Cost of holding a packet for one time slot\\
\hline
$C(Q_t,a_t)$ & Immediate cost if action $a_t$ is taken at time $t$ when the system is in state $Q_t$\\
\hline
$V(\theta)$ & Time average cost under the policy $\theta$ \\
\hline
$P_{a_t}(Q_t, Q_{t+1})$ & Transition probability from state $Q_t$ to $Q_{t+1}$ when action $a_t$ is chosen\\
\hline
$V_{\alpha}(i,j,\theta)$ & Total expected discounted cost under the policy $\theta$ when the initial state is $(i,j)$ \\
\hline
$V_{\alpha}(i,j)$ & Minimum total expected discounted cost when the initial state is $(i,j)$, i.e., $\min_{\theta} V_{\alpha}(i,j,\theta)$\\
\hline
$v_{\alpha}(i,j)$ & Difference of the minimum total expected discounted cost between the states $(i,j)$ and $(0,0)$, i.e., $V_{\alpha}(i,j)-V_{\alpha}(0,0)$\\
\hline
$V_{\alpha,n}(i,j)$ & Iterative definition for the optimality equation of $V_{\alpha}(i,j)$ \\
\hline
$\calV_{\alpha}(i,j,a)$ & $V_{\alpha}(i,j)= \min_{a \in \{0,1\}} \calV_{\alpha}(i,j,a)$, which is the optimality equation of $V_{\alpha}(i,j)$\\
\hline
$\Delta \calV(i,j)$& $\calV_{\alpha}(i,j,1)-\calV_{\alpha}(i,j,0)$ \\
\hline
\end{tabular}
\caption{Notation table}
\label{table:notation}
\end{table}

%% file: analysis.tex


\section{Should we maintain waiting time information?}

As described in the previous section, our goal is to obtain the  average-optimal policy.  To that end, we first find the space of possible policies and then identify the average-optimal policy within this space. Our first question is: what is the appropriate state space? Is it just queue length, or should we also consider waiting time?

Intuition tells us that if a packet has not been waiting for a long time then perhaps it could afford to wait a little more, but if a packet has waited for long, it might be better to just transmit it. That seems logical considering that we try our best to code but we cannot wait too long because it hurts in terms of holding costs. It is easy to keep track of waiting time information using time-stamps on packets when they are issued. Let $T^{(i)}$ be the arrival time of $i^{\text{th}}$ packet and $\mathcal{D}^{(i)}_{\theta}$ be its delay (i.e., the waiting time before it is transmitted) while policy $\theta$ is applied. We also denote by $\mathcal{T}_{t, \theta}$  the number of transmissions by time $t$ under policy $\theta$.
Then  Eq. (\ref{eq:long-time-average-cost}) can be written as
\bea
V(\theta)=\lim_{K \rightarrow \infty} \frac{1}{K+1} \expE_{\theta} \left[ \sum_{i: T^{(i)} \leq K} C_H \mathcal{D}^{(i)}_{\theta} + C_T  \mathcal{T}_{K,\theta}\right].
\eea
Would we be making better decisions by also keeping track of waiting times of each packet?   We answer this question in Proposition~\ref{ref:proposition2} that requires the following lemma, which indeed holds for generic MDPs \cite{MDP-Puterman}.


\begin{lemma}[\cite{MDP-Puterman}, Theorem 5.5.3]
For an MDP{$\{(Q_t, a_t), t \ge 0\}$}, given any randomized history dependent policy and starting state, there exists a randomized Markov policy with the same long-run average cost. \label{lemma:history-policy-imply-rand}
\end{lemma}


\begin{proposition}\label{ref:proposition2} \
\begin{itemize}
	\item[(i)] For the MDP$\{(Q_t, a_t), t \ge 0\}$, if there exists a randomized history dependent policy that is average-optimal then there exists a randomized Markov policy $\theta^* \in \Pi^{\text{MR}}$ that minimizes $V(\theta)$. 
	\item[(ii)] Further, one cannot find a policy which also uses waiting time information that would yield a better solution than $V(\theta^*)$.
\end{itemize}
\end{proposition}

\begin{proof}
The first result immediately follows from Lemma~\ref{lemma:history-policy-imply-rand}. Now, we focus on the second result. We notice that knowing the entire history of states (i.e., the number of packets in the queues) and actions one can always determine the history of waiting times as well as the current waiting times of all packets. Therefore the average-optimal policy $\theta^\prime$ that uses waiting time information is equivalent to a history dependent policy. From Lemma~\ref{lemma:history-policy-imply-rand}, we can always find a randomized Markov policy that yields the same average-optimal solution as $V(\theta^\prime)$.
\end{proof}

\section{Structure of the Average-Optimal Policy - Stationary and Deterministic Property}
In the previous section, we showed that there exists an  average-optimal policy that does not include the waiting time in the state of the system. Next, we focus on queue length based and randomized Markov policies, as well as determine the structure of the average-optimal policy. In this section, we will show that there exists an average-optimal policy that is stationary and deterministic. 

We begin by considering the infinite horizon $\alpha$-discounted cost case, where $0 < \alpha < 1,$ which we then tie to the average cost case.  This method \add[AS]{is typically used} in the MDP literature (e.g., \cite{stationary-policy-Sennott}), where the conditions for the structure of the average-optimal policy usually rely on the results of the infinite horizon $\alpha$-discounted cost case.  For our \mbox{MDP$\{(Q_t, a_t), t \ge 0\}$,} the total expected discounted cost incurred by a policy $\theta \in \Pi^{\text{HR}}$ is
\bea
V_{\alpha}(i,j,\theta) = \expE_\theta \left[ \sum_{t=0}^\infty \alpha^t C(Q_t, a_t) | Q_0 = (i,j) \right]. \label{eq:discount-cost}
\eea
In addition, we define $V_{\alpha}(i,j) = \min_{\theta} V_{\alpha}(i,j,\theta)$ as well as $v_{\alpha}(i,j)=V_{\alpha}(i,j)-V_{\alpha}(0,0)$. Define the $\alpha$-\textit{optimal policy} as the policy $\theta$ that minimizes \add[AS]{$V_{\alpha}(i,j,\theta)$}.\note[AS]{It seems that $\alpha$-optimal policy depends on the initial state?}

\subsection{Preliminary results}
In this subsection, we introduce the important properties of $V_{\alpha}(i,j)$, which  are mostly based on the literature \cite{stationary-policy-Sennott}.
We first show that $V_{\alpha}(i,j)$ is finite (Proposition~\ref{prop:finte-alpha-discount}) and then introduce the \textit{optimality equation} of $V_{\alpha}(i,j)$ (Lemma~\ref{proposition:optimality-eq}). 

\begin{proposition}
If \mbox{$\expE[\arrival_i] < \infty$} for $i=1,2$, then \mbox{$V_{\alpha}(i,j)< \infty$} for every state $(i,j)$ and $\alpha$.
\label{prop:finte-alpha-discount}
\end{proposition}

\begin{proof}
Let  $\tilde{\theta}$ be a stationary policy of  \change[AS]{always doing nothing}{waiting} (i.e., $a_t=0$ for all $t$) \annote[AS]{in each time slot}{even if in both  queues are not empty - not clear }. By definition of optimality, $V_{\alpha}(i,j) \leq V_{\alpha}(i,j,\tilde{\theta})$. Hence, if $V_{\alpha}(i,j,\tilde{\theta}) < \infty$, then $V_{\alpha}(i,j) < \infty$. Note that
\bean
V_{\alpha}(i,j,\tilde{\theta}) &=& \expE_{\tilde{\theta}}\Bigl[\sum_{t=0}^{\infty} \alpha^t C(Q_t, a_t) | Q_0=(i,j)\Bigr] \\
    &=&  \sum_{t=0}^{\infty} \alpha^t C_H \left(i+j +  t \expE [\arrival_1+\arrival_2]\right) \\
    &=& \frac{C_H (i+j)}{1-\alpha}+ \frac{\alpha C_H}{(1-\alpha)^2}  \expE[\arrival_1+\arrival_2] < \infty.
\eean
\end{proof}

The next lemma follows from Propositions 1 in \cite{stationary-policy-Sennott} and the fact that $V_{\alpha}(i,j)$ is finite (by Proposition~\ref{prop:finte-alpha-discount}). 

\begin{lemma}[\cite{stationary-policy-Sennott}, Proposition 1] \label{proposition:optimality-eq}
 \add[AS]{If \mbox{$\expE[\arrival_i] < \infty$} for $i=1,2$,  then the optimal expected discounted cost} $V_{\alpha}(i, j)$ satisfies the following \textit{optimality equation}:
\bea
V_\alpha(i,j) &=& \min_{a \in \{0,1\}} [ C_H ([i-a]^+ + [j-a]^+) + C_T a + \nonumber\\ 
&&\hspace{1cm}\alpha \sum_{k=0}^{\infty}\sum_{l=0}^{\infty} P_a\bigl((i,j),(k,l) \bigr) V_\alpha(k,l) ].  \label{eq:optimality-eq}
\eea
Moreover, the stationary policy that realizes the minimum of right hand side of  (\ref{eq:optimality-eq}) will be an $\alpha$-optimal policy.
\end{lemma}

\add[AS]{We define} $V_{\alpha, 0}(i, j)= 0$ and  for $n \ge 0$,
\begin{align}
V_{\alpha, n+1}(i,j) &= \min_{a \in \{0,1\}} [ C_H ([i-a]^+ + [j-a]^+) + C_T a + \nonumber\\
&\hspace{1cm}\alpha \sum_{k=0}^{\infty}\sum_{l=0}^{\infty} P_a\bigl((i,j),(k,l)\bigr) V_{\alpha, n}(k,l) ]. \label{eq:optimality-itr}
\end{align}
Lemma~\ref{proposition:optimality-eq1} below follows from Proposition 3 in \cite{stationary-policy-Sennott}.
 
\begin{lemma}[\cite{stationary-policy-Sennott}, Proposition 3] \label{proposition:optimality-eq1}
 $V_{\alpha, n}(i, j) \rightarrow V_{\alpha}(i, j)$ as $n \rightarrow \infty$  for every $i$, $j$, and $\alpha$.
\end{lemma}

\change[AS]{We remark that}{ Eq.} (\ref{eq:optimality-itr}) will be helpful for identifying the properties of $V_{\alpha}(i, j)$, e.g., to prove that  $V_{\alpha}(i,j)$ is a non-decreasing function. \remove[AS]{as in the following Lemma.} 

\begin{lemma} \label{lemma:increasing function}
$V_{\alpha}(i, j)$ is a non-decreasing function with respect to (w.r.t.) $i$ for fixed $j$, and vice versa.
\end{lemma}

\begin{proof}
The proof is by induction on $n$ in Eq. (\ref{eq:optimality-itr}). The result clearly holds for  $V_{\alpha, 0}(i, j)$. Now, assume that $V_{\alpha, n}(i, j)$
is non-decreasing. First, note that $C_H ([i-a]^+ + [j-a]^+)+ C_T a$ is a non-decreasing function of $i$ and $j$ (since $C_H$ is non-negative). Next, we note that 
\begin{align}
& \alpha \sum_{k=0}^{\infty}\sum_{l=0}^{\infty} P_a\bigl((i,j),(k,l)\bigr) V_{\alpha, n}(k,l)  \nonumber\\
=& \alpha \sum_{r=0}^{\infty}\sum_{s=0}^{\infty} p^{(1)}_r p^{(2)}_s V_{\alpha, n}( [i-a]^++r, [j-a]^++s),\nonumber
\end{align}
which is also a non-decreasing function in $i$ and $j$ separately due to the inductive assumption. Since the sum and the minimum (in Eq. (\ref{eq:optimality-itr})) of non-decreasing functions are a non-decreasing function, we conclude that $V_{\alpha, n+1}(i, j)$ is a non-decreasing function as well. 
\end{proof}

\add[AS]{The next two lemmas, which can be proven via the similar arguments in \cite{stationary-policy-Sennott}, specify the conditions for the existence of the optimal stationary and deterministic policy.}


\begin{lemma}[\cite{stationary-policy-Sennott}, Theorem (i)] \note[AS]{This lemma is important - can we change it to a theorem?}
\label{lemma:stationary-policy-exist} 
There exists a stationary and deterministic policy that is average-optimal for the MDP$\{(Q_t,a_t), t \ge 0\}$ if the following conditions are satisfied: 
\begin{itemize}
	\item[(i)] $V_{\alpha}(i,j)$ is finite for all $i$, $j$, and $\alpha$;\note[AS]{should it be finite for a given $\alpha$ or for any value of $\alpha$?}
	\item[(ii)] There exists a nonnegative $N$ such that $v_{\alpha}(i,j) \geq -N$ for all $i$, $j$, and $\alpha$; 
	\item[(iii)] There exists a nonnegative $M_{i,j}$ such that $v_{\alpha}(i,j) \leq M_{i,j}$ for every $i$, $j$, and $\alpha$. Moreover, for each state $(i,j)$ there is an action $a(i,j)$\note[AS]{we might have some definition for the action before - and this seems to be different - should we introduce a new one? Also, $G$ might not be a good notation, since it is usually used for graphs.} such that $\sum_{k=0}^{\infty}\sum_{l=0}^{\infty} P_{a(i,j)}\bigl((i,j),(k,l)\bigr)M_{k,l} < \infty$.\note[AS]{should we specify the range of $j$ and $k$?}
\end{itemize}
\end{lemma}

\begin{lemma}[\cite{stationary-policy-Sennott}, Proposition 5]
\label{lemma:stationary-policy-exist2}
Assume there exists a stationary policy $\theta$ inducing an irreducible and ergodic Markov chain with the following properties: there exists a nonnegative function $F(i,j)$ and a finite nonempty subset $G \subseteq (\mathbb{N}\cup \{0\})^2$ such that  for $(i,j) \in (\mathbb{N}\cup \{0\})^2 - G$ it holds that
\begin{align}
\sum_{k=0}^{\infty} \sum_{l=0}^{\infty} P_{a(\theta)}((i,j),(k,l))F(k,l)-F(i,j)  \leq -C((i, j), a(\theta)), \label{eq:average-cost-opt-stationary}
\end{align}
where $a(\theta)$ is the action when the policy $\theta$ is applied. Moreover, for $(i,j) \in G$ it holds that 
$$\sum_{k=0}^{\infty}\sum_{l=0}^{\infty} P_{a(\theta)}((i,j), (k,l)) F(k,l) < \infty.$$
Then, the condition (iii) in Lemma \ref{lemma:stationary-policy-exist} holds.
\end{lemma}

\subsection{Main result}
Using lemmas \ref{lemma:stationary-policy-exist} and \ref{lemma:stationary-policy-exist2}, we show next that the MDP defined in this paper has an average-optimal policy that is stationary and deterministic.
\begin{theorem} \label{thm:stationary}
For the MDP$\{(Q_t, a_t), t \ge 0\}$, there exists a stationary and deterministic policy $\theta^*$ that minimizes $V(\theta)$ if $\expE[\arrival_i^2] < \infty$ and $\expE[\arrival_i] < 1$ for  $i=1,2$.  
\end{theorem}

\begin{proof}
As described earlier it is sufficient to show that the three conditions in Lemma~\ref{lemma:stationary-policy-exist} are satisfied. Proposition~\ref{prop:finte-alpha-discount} implies that the condition (i) holds, while the condition (ii) is satisfied due to  Lemma~\ref{lemma:increasing function} (i.e., $N=0$ in Lemma \ref{lemma:stationary-policy-exist}).  We denote by $\tilde{\theta}$ the stationary policy of transmitting \add[AS]{at each time slot}.\note[AS]{moved the definition.} We use this policy for each of the three cases described blow and show that condition (iii) of Lemma~\ref{lemma:stationary-policy-exist} holds.\note[AS]{Are you going to show that the condition of Lemma~8 holds?}

\textbf{Case (i):}  $p^{(i)}_0+p^{(i)}_1<1$ for $i=1,2$, i.e., the probability that two or more packets arrive for each time slot is non-zero.
This policy $\tilde{\theta}$ results in an irreducible and ergodic Markov chain, and therefore Lemma \ref{lemma:stationary-policy-exist2} can be applied. Let $F(i, j)= B (i^2+ j^2)$ for some positive $B$. Then, for all states $(i,j) \in (\mathbb{N} \cup \{0\})^2 -\{(0,0),(0,1), (1,0)\}$, it holds that
\begin{align*}
&\sum_{k=0}^{\infty}\sum_{l=0}^{\infty} P_{a(\tilde{\theta})}\left((i,j), (k,l)\right)\left[F(k,l)-F(i,j)\right] \\
=& \sum^{\infty}_{r=0} \sum^{\infty}_{s=0}P_1\left((i,j), ([i-1]^++r,[ j-1]^++s)\right) \cdot\\
&\hspace{1.2cm} \left[ F([i-1]^++r, [j-1]^++s) - F(i,j)\right] \\
=& \sum^{\infty}_{r=0} \sum^{\infty}_{s=0} p^{(1)}_r p^{(2)}_s B\bigl[ 2i(r-1)+ (r-1)^2+\\ 
& \hspace{1.2cm}2j(s-1)+ (s-1)^2\bigr] \\
=& 2B\Bigl( i ( \expE[\arrival_1]  - 1) + j( \expE[\arrival_2]  - 1)\Bigr)+ \\
&B  \Bigl( \expE[(\arrival_1-1)^2] +   \expE[(\arrival_2-1)^2]\Bigr). 
\end{align*}

Note that $\expE[\arrival_i] < 1$, hence $2B( \expE[\arrival_i]-1) < -C_H$ for sufficiently large $B$. Moreover, since $\expE[\arrival_i^2] < \infty$, it holds that 
\begin{eqnarray*}
&&\sum_{k=0}^{\infty}\sum_{l=0}^{\infty} P_{a(\tilde{\theta})} \left((i,j), (k,l)\right)\left[F(k,l)-F(i,j)\right]  \\
&\leq& -C((i,j), a(\tilde{\theta})),
\end{eqnarray*}
when $i, j$ are large enough, where $$C((i,j), a(\tilde{\theta})) = C_H([i-1]^++[j-1]^+)+C_T.$$ 

We observe that there exists a finite set $G$  that contains states $\{(0,0), (0,1), (1,0)\}$ such that Eq. (\ref{eq:average-cost-opt-stationary}) is satisfied for $(i,j) \in (\mathbb{N} \cup \{0\})^2 - G$.
 Then, for $(i,j) \in G$, it holds that
\begin{align*}
&\sum_{k=0}^{\infty}\sum_{l=0}^{\infty} P_{a(\tilde{\theta})}\left((i,j), (k,l)\right) F(k,l) \\
=& B \sum^{\infty}_{r=0} \sum^{\infty}_{s=0} p^{(1)}_r p^{(2)}_s \left[ ([i-1]^++r)^2 + ([j-1]^++s)^2\right] \\
=& B \Bigl\{(i-1)^2+ 2[i-1]^+ \expE[\arrival_1] + \expE[\arrival_1^2] + \\
&(j-1)^2+ 2[j-1]^+ \expE[\arrival_2] + \expE[\arrival_2^2] \Bigr\}< \infty.
\end{align*}
Therefore, the condition of Lemma~\ref{lemma:stationary-policy-exist2} is satisfied, which implies, in turn, that  condition  (iii) in Lemma~\ref{lemma:stationary-policy-exist} is satisfied as well. 

%

\begin{figure}
\centering
\includegraphics[width=0.38 \textwidth]{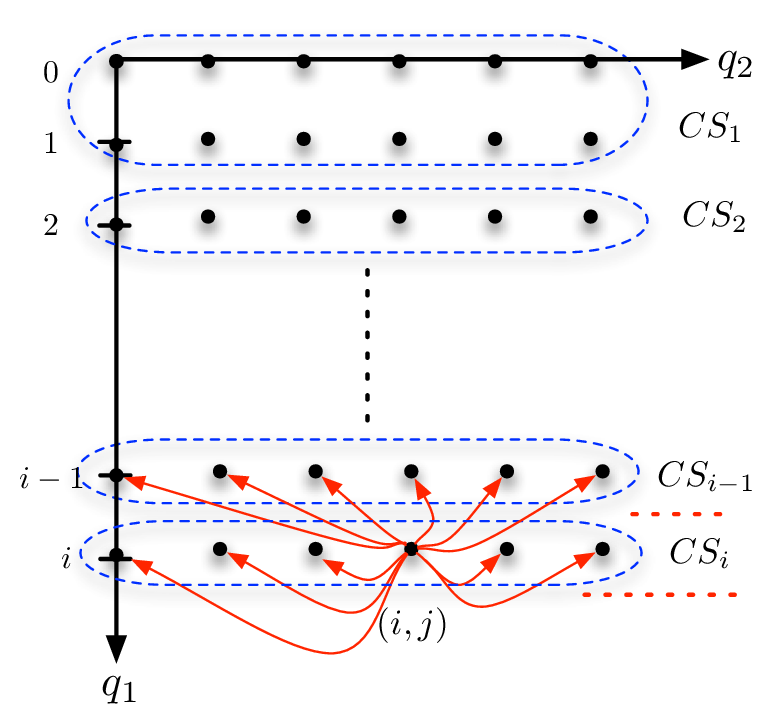}
\caption{Case (ii) in the proof of Theorem \ref{thm:stationary}: state $(i,j)$ can only transit to the states in the $CS_i$ and $CS_{i-1}$.}
\label{fig:comm-class}
\end{figure}

\textbf{Case (ii):} $p^{(1)}_0+p^{(1)}_1=1$ and $p^{(2)}_0+p^{(2)}_1 < 1$. Note that $\tilde{\theta}$  \annote[AS]{results in a reducible Markov chain}{note the change in language}. That is,  \add[AS]{there are several communicating classes} as depicted in Fig. \ref{fig:comm-class}. 

We define the classes $CS_1= \{(a,b): a=0, 1 \text{\,\,and\,\,} b\in \mathbb{N} \cup \{0\} \}$ and  $CS_i=\{(a,b): a=i, b\in \mathbb{N} \cup \{0\} \}$ for $i \geq 2$. Then each $CS_i$ is a  communicating class under policy $\tilde{\theta}$. The states in $CS_1$ are positive-recurrent, and  each $CS_i$ for $i \geq 2$ is a transient class. For $i \geq 2$, let $\overline{C}_{i,j}$ be the expected cost of  the passage from state $(i,j)$ (in class $CS_i$) to the next class $CS_{i-1}$. Note that state $(i,j)$ has the probability of $p^{(1)}_0$ to escape to class $CS_{i-1}$ and $p^{(1)}_1$ to remain in class $CS_{i}$. Now $C((i, Q^{(2)}_t),1) = C_T + C_H( [i-1]^+ + [Q^{(2)}_t-1]^+)$. By considering all the possible paths to escape from state $(i,j)$, we can compute $\overline{C}_{i,j}$ as follows:
\begin{align*}
\overline{C}_{i,j}=& \expE \left[ \sum_{k=0}^{\infty} (p^{(1)}_1)^k p^{(1)}_0  \sum^{k}_{t=0} C((i, Q^{(2)}_t),1)  | (Q^{(1)}_0,Q^{(2)}_0) =(i, j) \right] \\ 
=& p^{(1)}_0  \expE \left[  \sum^{\infty}_{t=0} C((i, Q^{(2)}_t),1) \sum^{\infty}_{k=t}  (p^{(1)}_1)^k  | (Q^{(1)}_0,Q^{(2)}_0) =(i, j)\right] \\
=&\expE \left[\sum^{\infty}_{t=0} (p^{(1)}_1)^t  C((i, Q^{(2)}_t),1) | (Q^{(1)}_0,Q^{(2)}_0) =(i, j) \right]. 
\end{align*}
We observe that $\overline{C}_{i,j}$ can be viewed as the total expected $p^{(1)}_1$-discounted cost of the system.
Following the arguments similar to these in the proof of Proposition~\ref{prop:finte-alpha-discount}, we conclude that $\overline{C}_{i,j} < \infty$.

We denote the expected cost of a first passage from state $(i,j)$ to $(k,l)$ by  $\overline{C}_{(i,j),(k,l)}.$  Proposition~4 in \cite{stationary-policy-Sennott} implies  that $\overline{C}_{(1,j),(0,0)} < \infty$ for any $j$, where the intuition is that the expected traveling time from state $(1,j)$ to $(0,0)$ is finite due to the positive recurrence of $CS_{1}$.
Let $T_0=\min \{t \geq 1: (Q^{(1)}_t,Q^{(2)}_t)=(0,0)\}$ and for $i \geq 1$, $T_i =\min \{t \geq 1: Q^{(1)}_t=i\}$ with the corresponding state $(Q^{(1)}_{T_i}, Q^{(2)}_{T_i}) = (i, \tilde{j}_i)$.\note[AS]{corresponding to what? What is $\tilde{j}_i$?} Since
\bean
\overline{C}_{(i,j),(0,0)}= \overline{C}_{i,j} + \sum_{k=1}^{i-2} \overline{C}_{i-k, \tilde{j}_{i-k}} + \overline{C}_{(1,\tilde{j}_1), (0,0)},
\eean
we conclude that  $\overline{C}_{(i,j),(0,0)} < \infty$.

Let $\hat{\theta}$, be a policy that always transmits until time slot $T_0$ after which  the $\alpha$-optimal policy is employed. Then, $V_{\alpha}(i,j)$ can be bounded by
\bean
V_{\alpha}(i,j) \leq && \expE_{\hat{\theta}} \left[ \sum_{t=0}^{T_{0}-1} \alpha^t C(Q_t, a_t) | Q_0 = (i,j) \right] +  \\
&&\expE_{\hat{\theta}}  \left[ \sum_{t=T_0}^{\infty} \alpha^t C(Q_t, a_t) | Q_0 = (i,j) \right] \\
\leq&& \overline{C}_{(i,j),(0,0)} + V_{\alpha}(0,0).
\eean
We show that  condition (iii) of  Lemma~\ref{lemma:stationary-policy-exist} is satisfied by choosing $M_{i,j}=\overline{C}_{(i,j),(0,0)}$. In particular, it holds that $v_{\alpha}(i,j)=V_{\alpha}(i,j)-V_{\alpha}(0,0) \leq M_{i,j}$ and $M_{i,j} < \infty$. Moreover, $\sum_{k=0}^{\infty}\sum_{l=0}^{\infty} P_1\bigl( (i,j),(k,l) \bigr) M_{k,l} = \sum_{k=0}^{\infty}\sum_{l=0}^{\infty} P_1\bigl((i,j),(k,l) \bigr)\overline{C}_{(k,l),(0,0)} \leq \overline{C}_{(i,j),(0,0)} < \infty$. 

\textbf{Case (iii):}  $p^{(i)}_0+p^{(i)}_1=1$ for $i=1,2$, i.e., Bernoulli arrivals to both queues.  Note that in this case $\tilde{\theta}$  also \annote[AS]{results in a reducible Markov chain.}{Also not clear to me.} The proof is similar  to case (ii) - we define $M_{i,j}=\overline{C}_{(i,j),(0,0)}$, and show that $\overline{C}_{(i,j),(0,0)}$ is finite for this case. 
\end{proof}

According to Borkar \cite{borkar-convex}, it is possible to find the randomized policy that is closed to the average-optimal by applying linear programming methods for an MDP of a very generic setting, where randomized stationary  policies are average-optimal.  However, since the average-optimal policy has further been shown in Theorem \ref{thm:stationary} to be deterministic, in the next section we  investigate the  structural properties of the average-optimal policy and using a Markov-chain based enumeration to find the average-optimal polity that would be deterministic stationary.

\section{Structure of the Average-Optimal Policy - Threshold Based}
Now that we know the average-optimal policy is stationary and deterministic, the question is how do we find it?  If we know that the average-optimal policy satisfies the structural properties, then it is possible to search through the space of stationary deterministic policies and obtain the optimal one. We will study the $\alpha$-optimal policy first and then discuss how to correlate it with the average-optimal policy. Before investigating the general i.i.d. arrival model, we study a special case, namely Bernoulli process.  Our objective is to determine the $\alpha$-optimal policy for the Bernoulli arrival process.

\begin{lemma} \label{thm:bernoulli-threshold}
For the  i.i.d. Bernoulli arrival process and the system starting from the empty queues, the $\alpha$-optimal policy is of threshold type. In particular, there exist  optimal thresholds $L^*_{\alpha, 1}$ and $L^*_{\alpha,2}$ so that the optimal deterministic action in state $(i,0)$ is to wait if $i \leq L^*_{\alpha,1}$, and to transmit without coding if $i > L^*_{\alpha,1}$; while in state $(0,j)$ is to wait if $j \leq L^*_{\alpha,2}$, and to transmit without coding if $j > L^*_{\alpha,2}$. 
\end{lemma}
\begin{proof}
\add[AS]{We define}  
\begin{eqnarray*}
&&\calV_{\alpha}(i,0,a) \\
&=&  C_H ([i-a]^+) + C_T a + \alpha \sum_{k,l} P_a\bigl((i,0),(k,l)\bigr) V_\alpha(k,l).
\end{eqnarray*}
Then, $$V_{\alpha}(i,0) = \min_{a \in \{0,1\}} \calV_{\alpha}(i,0,a).$$ Let $L^*_{\alpha,1} = \min \{i \in \mathbb{N} \cup \{0\}: \calV_{\alpha}(i,0, 1) > \calV_{\alpha}(i, 0,0) \}-1$. Then the optimal stationary and deterministic action (for the total expected $\alpha$-discounted cost) is $a_t=0$ for the states $(i,0)$ with $i \leq L^*_{\alpha,1}$, and $a_t=1$ for the state $(L^*_{\alpha,1}+1,0)$.\note[AS]{Should we explain why? It is not immediately clear from the discussion.} Note that we do not need to define the policy of states $(i,0)$ for  $i>L^*_{\alpha,1}+1$, since they are not accessible as $(L^*_{\alpha,1}+1,0)$ only transits to $(L^*_{\alpha,1},0)$, $(L^*_{\alpha,1}+1,0)$, $(L^*_{\alpha,1},1)$, and $(L^*_{\alpha,1}+1,1)$. The similar argument is applicable for the states $(0, j)$. Consequently, there exists a  policy of threshold type that is $\alpha$-optimal.
\end{proof}

\subsection{General i.i.d. arrival process}
For the i.i.d. Bernoulli arrival process, we have just shown that the  $\alpha$-optimal policy is threshold based. Our next objective is to extend \change[AS]{the ideas}{this result} to any i.i.d. arrival process. 
We define that $\calV_{\alpha}(i,j,a) = C_H \left([i-a]^+ + [j-a]^+\right) + C_T \cdot a + \alpha \expE  [V_{\alpha} \left([i-a]^+ +\arrival_1, [j-a]^+ + \arrival_2 \right)] $. Moreover, let $\calV_{\alpha,n}(i,j,a) = C_H \left([i-a]^+ + [j-a]^+\right) + C_T a + \alpha \expE  [V_{\alpha,n} \left([i-a]^+ +\arrival_1, [j-a]^+ + \arrival_2 \right)] $. Then Eq. (\ref{eq:optimality-eq}) can be written as $V_{\alpha}(i,j)=\min_{a \in \{0,1\}} \calV_{\alpha}(i,j, a)$, while Eq. (\ref{eq:optimality-itr})  can be written as $V_{\alpha,n+1}(i,j)=\min_{a \in \{0,1\}} \calV_{\alpha, n}(i,j, a)$. For every discount factor $\alpha$, we want  to show that there exists an $\alpha$-optimal policy that is of threshold type. To be precise, let the $\alpha$-optimal policy for the first dimension be $a^*_{\alpha, i}=\min{\{a' \in \argmin_{a \in \{0,1\}} \calV_{\alpha}(i,0,a)}\}$,\footnote{This notation also used in \cite{MDP-Puterman} combines two operations: First we let $\Lambda=\{a \in\{0,1\}: \min \calV_{\alpha,n}(i,0,a)\}$, and then do $\min \Lambda$. In other words, we choose $a=0$ when both $a=0$ and $a=1$ result in the same $\calV_{\alpha,n}(i,j,a)$.} and we will show that $a^*_{\alpha, i}$ is non-decreasing as $i$ increases, and so is the second dimension. We start with a number of definitions that describe the properties of $V_{\alpha}(i,j)$.


\begin{define}[\cite{monotone-policy-survey}, Submodularity]
A function $f: (\mathbb{N} \cup \{0\})^2 \rightarrow \mathbb{R}$ is submodular if for all $i, j \in \mathbb{N} \cup \{0\}$
\bean
f(i,j)+ f(i+1, j+1) \leq f(i+1,j)+ f(i, j+1).
\eean
\end{define}

\begin{define}[$\calk$-Convexity] \label{def:convex}
A function $f: (\mathbb{N} \cup \{0\})^2 \rightarrow \mathbb{R}$ is $\calk$-convex (where $\mathcal{K} \in \mathbb{N}$) if for every $i, j \in \mathbb{N} \cup \{0\}$
\bean
f(i+\calk,j) - f(i,j ) \leq f(i+\calk+1,j) - f(i+1, j); \\
f(i,j+\calk) - f(i,j ) \leq f(i, j+\calk+1) - f(i, j+1).
\eean
\end{define}

\begin{define}[$\calk$-Subconvexity] \label{def:subconvex}
A function $f: (\mathbb{N} \cup \{0\})^2\rightarrow \mathbb{R}$ is $\calk$-subconvex (where $\mathcal{K} \in \mathbb{N}$) if for all $i, j \in \mathbb{N} \cup \{0\}$
\bean
f(i+\calk,j+\calk) - f(i,j ) \leq f(i+\calk+1, j+\calk) - f(i+1, j); \\
f(i+\calk,j+\calk) - f(i,j ) \leq f(i+\calk, j+\calk+1) - f(i, j+1).
\eean
\end{define}

\begin{remark}
If a function $f: (\mathbb{N} \cup \{0\})^2 \rightarrow \mathbb{R}$ is submodular and $\calk$-subconvex, then it is $\calk$-convex, and for every $r \in \mathbb{N} $ with $1 \leq r<\calk$,
\bean
f(i+ \calk, j+r)-f(i,j) \leq f(i+\calk+1,j+r) - f(i+1,j); \\
f(i+r , j+\calk)-f(i,j) \leq f(i+r,j+\calk+1) - f(i,j+1). 
\eean
\end{remark}

For simplicity, we will ignore $\calk$ in definitions \ref{def:convex} and \ref{def:subconvex} when $\calk=1$. We will show in Subsection \ref{subsection:main-result} that $V_{\alpha}(i,j)$ is non-decreasing, submodular, and subconvex, that result in the threshold base of $\alpha$-optimal policy. Note that the definition of $\calk$-Convexity (Definition~\ref{def:convex}) is dimension-wise, which is different from the definition of convexity for the continuous function in two dimensions.   

\subsection{Proof overview}
Before the technical proofs  in Subsection \ref{subsection:main-result}, in this subsection, we overview why submodularity and subconvexity  of $V_{\alpha}(i,j)$ lead to the $\alpha$-optimality of the threshold based policy. 

\begin{itemize}
	\item \textit{To show that $\alpha$-optimal policy is monotonic w.r.t. state $(i,0)$, it suffices to show that $\calV_{\alpha}(i,0,1)-\calV_{\alpha}(i,0,0)$ is a non-increasing function w.r.t. $i$:}  Suppose it is true that $\calV_{\alpha}(i+1,0,1)-\calV_{\alpha}(i+1,0,0) \leq \calV_{\alpha}(i,0,1)-\calV_{\alpha}(i,0,0)$. We observe that if the $\alpha$-optimal policy for state $(i,0)$ is $a^*_{\alpha, i}=1$, i.e., $\calV_{\alpha}(i,0,1) - \calV_{\alpha}(i,0,0) \leq 0$,  then the $\alpha$-optimal policy for state $(i+1,0)$ is also $a^*_{\alpha, i+1}=1$. Similarly, if the $\alpha$-optimal policy for state $(i+1,0)$ is $a^*_{\alpha, i+1}=0$ then the $\alpha$-optimal policy for state $(i,0)$ is $a^*_{\alpha,i}=0$. 
	
	\item  \textit{In oder to prove that $\calV_{\alpha}(i,0,1)-\calV_{\alpha}(i,0,0)$ is non-increasing, it is sufficient to show that $V_{\alpha}(i,j)$ is convex:}  When $i \geq 1$, the claim is true since
\begin{align*}
&\calV_{\alpha}(i,0,1)-\calV_{\alpha}(i,0,0)\\
=& C_T- C_H+ \alpha \expE[V_{\alpha}(i-1+\arrival_1,\arrival_2)- V_{\alpha}(i+\arrival_1,\arrival_2)]. 
\end{align*}

	\item \textit{Similarly, to show that $\alpha$-optimal policy of state $(i,j)$ is monotonic w.r.t.  $i$ for  fixed $j$ and vice versa,  it suffices to show that $V_{\alpha}(i,j)$ is subconvex:} When $i,j \geq 1$, we observe that
\begin{align*}
&\calV_{\alpha}(i,j,1)-\calV_{\alpha}(i,j,0) \\
=& C_t- 2C_h+ \alpha \expE[V_{\alpha}(i-1+\arrival_1,j-1+\arrival_2)- \\
&\hspace{2.3cm}V_{\alpha}(i+\arrival_1,j+\arrival_2)]. 
\end{align*}

	\item \textit{To show $V_{\alpha}(i,j)$ is convex and subconvex, we need $V_{\alpha}(i,j)$ is submodular:} We intend to prove the convexity and subconvexity of $V_{\alpha}(i,j)$ by induction, which will require the relation between $V_{\alpha}(i,j)+V_{\alpha}(i+1,j+1)$ and $V_{\alpha}(i+1,j)+V_{\alpha}(i,j+1)$. There will be two choices: (i) $V_{\alpha}(i,j)+V_{\alpha}(i+1,j+1) \leq V_{\alpha}(i+1,j)+ V_{\alpha}(i,j+1)$, or (ii) $V_{\alpha}(i,j)+V_{\alpha}(i+1,j+1) \geq V_{\alpha}(i+1,j)+ V_{\alpha}(i,j+1)$.  We might assume that $V_{\alpha}(i,j)$  satisfies (i). Then (i) and the subconvexity of $V_{\alpha}(i,j)$  implies the convexity of $V_{\alpha}(i,j)$. In the contrary,  the convexity of $V_{\alpha}(i,j)$ and (ii) lead to the subconvexity of $V_{\alpha}(i,j)$. In other words, both  choices are possible since they do not violate the convexity and subconvexity of $V_{\alpha}(i,j)$. Now we are going to argue that the choice (ii) is wrong. Suppose  the actions of $\alpha$-optimal policy for the states $(i,j)$, $(i+1,j)$, $(i,j+1)$,  $(i+1,j+1)$ are  $0,0,1,1$ respectively. If the choice (ii) is true, then when $i \geq 1$, we have
\begin{align*}
&C_H(i+j) +  \expE[V_{\alpha,n}(i+\arrival_1,j+\arrival_2)]+\\
& C_T+ C_H(i+j)+\expE[V_{\alpha,n}(i+\arrival_1, j+\arrival_2)]  \\
\geq& C_H(i+1+j)+ \expE[V_{\alpha,n}(i+1+\arrival_1, j+\arrival_2)]+\\
& C_T+C_H(i-1+j)+\expE[V_{\alpha,n}(i-1+\arrival_1,j+\arrival_2)].
\end{align*}
By simplifying the above inequality, we can observe  the contradiction to the fact that $V_{\alpha,n}(i,j)$ is convex. Therefore,  $V_{\alpha}(i,j)$ is submodular.

\end{itemize}

So far, we know that if we show $V_{\alpha}(i,j)$ is submodular and subconvex, then the $\alpha$-optimal policy of state $(i,j)$ is  non-decreasing separately in the direction of $i$ and $j$ (i.e., threshold type). Next, we briefly discuss how Lemmas \ref{lemma:subadditivity}-\ref{lemma:subconvex} and Theorem \ref{thm:thresold-discounted} in the next subsection work together. 
Theorem  \ref{thm:thresold-discounted} states that the $\alpha$-optimal policy is of threshold type, with the proof of induction on $n$ in Eq. (\ref{eq:optimality-itr}). First, we observe that $V_{\alpha,0}(i,j)$ is non-decreasing, submodular, and subconvex. Second, based on Lemma \ref{lemma:subadditivity} and  Corollary \ref{cor:itr-monotone-policy},  $\min\{a' \in \argmin_{a \in \{0,1\}} \calV_{\alpha,0}(i,j,a)\}$ is non-decreasing w.r.t. $i$ for fixed $j$, and vice versa. Third, according to Lemmas \ref{lemma:increasing function}, \ref{lemma:submodular}, and \ref{lemma:subconvex}, we know that $V_{\alpha,1}(i,j)$ is non-decreasing, submodular, and subconvex. Therefore, as $n$ goes to infinity, we  conclude that $V_{\alpha}(i,j)$ is non-decreasing, submodular, and subconvex, as well as $\min\{a' \in \argmin_{a \in \{0,1\}} \calV_{\alpha}(i,j,a)\}$ is non-decreasing w.r.t. $i$ for fixed $j$, and vice versa.

\subsection{Main results and proofs} \label{subsection:main-result}
\begin{lemma} \label{lemma:subadditivity}
Given $0 < \alpha <1$ and $n \in \mathbb{N} \cup \{0\}$. If $V_{\alpha,n}(i,j)$ is non-decreasing, submodular, and subconvex, then $\calV_{\alpha,n}(i,j,a)$ is submodular for $i$ and $a$ when $j$ is fixed, and so is for $j$ and $a$ when $i$ is fixed.
\end{lemma}

\begin{proof}
We define $\Delta \calV_{\alpha,n}(i,j)= \calV_{\alpha,n}(i,j,1)-\calV_{\alpha,n}(i,j,0)$. We claim that $\Delta \calV_{\alpha, n}(i,j)$ is non-increasing, i.e., $\Delta \calV_{\alpha,n}(i,j)$ is a non-increasing function w.r.t. $i$ while $j$ is fixed, and vice versa (we will focus on the \change{prior}{former} part). Notice that 
\bean
\Delta \calV_{\alpha,n}(i,j) &=& C_H([i-1]^+ + (j-1)^+)+ C_T+ \\
&&\alpha \expE[ V_{\alpha,n}([i-1]^+ + \arrival_1, [j-1]^++ \arrival_2 )]- \nonumber\\
        && C_H(i+j) - \alpha \expE[ V_{\alpha,n}(i+\arrival_1,j+\arrival_2)].
\eean
To be precise, when $i \geq 1$,
\begin{align}
\Delta \calV_{\alpha,n}(i,j) =&  C_T- 2C_H+ \alpha \expE[V_{\alpha,n}(i-1+\arrival_1,j-1+\arrival_2)- \nonumber\\
& V_{\alpha,n}(i+\arrival_1,j+\arrival_2)]  \,\,\,\text{for}\,\, j \geq 1; \label{eq:subadditivity1} \\
\Delta \calV_{\alpha,n}(i,j) =& C_T- C_H+ \alpha \expE[V_{\alpha,n}(i-1+\arrival_1,\arrival_2)- \nonumber\\
&V_{\alpha,n}(i+\arrival_1,\arrival_2)]  \,\,\, \text{for}\,\, j =0.\label{eq:subadditivity2}
\end{align}
Because of the subconvexity of $V_{\alpha,n }(i,j)$ in Eq. (\ref{eq:subadditivity1}), when $i \geq 1$ and $j \geq 1$, $\Delta \calV_{\alpha, n}(i,j)$ does not increase as $i$ increases. The same is for $i \geq 1$ and $j=0$ in Eq. (\ref{eq:subadditivity2}) due to the convexity of $V_{\alpha,n }(i,j)$.

\add{We proceed to establish} the boundary conditions. When $j \geq 1$,
\bean
\Delta \calV_{\alpha,n}(1,j) &=& C_T- 2C_H + \alpha \expE[V_{\alpha,n}(\arrival_1, j-1+\arrival_2)- \\
&&\hspace{2.5cm} V_{\alpha,n}(1+\arrival_1, j+\arrival_2)];  \\
\Delta \calV_{\alpha,n}(0,j) &=&  C_T- C_H+ \alpha \expE[V_{\alpha,n}(\arrival_1,j-1+\arrival_2)- \\
&&\hspace{2.5cm}V_{\alpha,n}(\arrival_1,j+\arrival_2)].
\eean
Note that $\expE[V_{\alpha,n}(1+ \arrival_1, j+\arrival_2)] \geq \expE[V_{\alpha,n}(\arrival_1,j+\arrival_2)]$ according to non-decreasing $V_{\alpha,n}(i,j)$ and then $\Delta \calV_{\alpha,n}(1,j) \leq \Delta \calV_{\alpha,n}(0,j)$ when $j \geq 1$. Finally, when $j=0$ we have
\bean
\Delta \calV_{\alpha,n}(1,0)&=& C_T- C_H + \alpha \expE[V_{\alpha,n}(\arrival_1, \arrival_2)-\\
&& \hspace{2.5cm}V_{\alpha,n}(1+\arrival_1, \arrival_2)]; \\
\Delta \calV_{\alpha,n}(0,0) &=&C_T.
\eean
Here, $\Delta \calV_{\alpha,n}(1,0) \leq \Delta \calV_{\alpha,n}(0,0)$ since $\expE[V_{\alpha,n}(\arrival_1,\arrival_2)- V_{\alpha,n}(1+\arrival_1,\arrival_2)] \leq 0$ as $V_{\alpha,n}(i,j)$ is non-decreasing. Consequently, $\Delta \calV_{\alpha,n}(i,j)$ is a non-increasing function w.r.t. $i$ while $j$ is fixed. 
\end{proof}

Submodularity of $\calV_{\alpha,n}(i,j,a)$ implies the monotonicity of the optimal minimizing policy \cite[Lemma 4.7.1]{MDP-Puterman} as described in the following Corollary. This property will simplify the proofs of Lemmas \ref{lemma:submodular} and \ref{lemma:subconvex}.
\begin{corollary} \label{cor:itr-monotone-policy}
Given $0 < \alpha <1$ and $n \in \mathbb{N} \cup \{0\}$. If $V_{\alpha,n}(i,j)$ is non-decreasing, submodular, and subconvex, then $\min\{a' \in \argmin_{a \in \{0,1\}} \calV_{\alpha,n}(i,j,a)\}$ is non-decreasing w.r.t. $i$ for fixed $j$, and vice versa.
\end{corollary}

\begin{lemma} \label{lemma:submodular}
Given $0 < \alpha <1$ and $n \in \mathbb{N} \cup \{0\}$. If $V_{\alpha,n}(i,j)$ is non-decreasing, submodular, and subconvex, then $V_{\alpha, n+1}(i,j)$ is submodular.
\end{lemma}
\begin{proof}
We intend to show that $V_{\alpha, n+1}(i+1,j+1)-V_{\alpha, n+1}(i+1,j) \leq V_{\alpha, n+1}(i,j+1)- V_{\alpha, n+1}(i,j)$ for all $i, j \in \mathbb{N} \cup \{0\}$. According to Corollary \ref{cor:itr-monotone-policy}, only 6 cases of $(a^*_{i,j}, a^*_{i+1,j}, a^*_{i,j+1}, a^*_{i+1,j+1})$ are considered, where $a^*_{i,j}= \min\{a' \in \argmin_{a \in \{0,1\}} \calV_ {\alpha, n}(i,j,a) \}$.

\textbf{Case (i):} if $(a^*_{i,j}, a^*_{i+1,j}, a^*_{i,j+1}, a^*_{i+1,j+1})=(1,1,1,1)$,  we claim that
\bean
&&\expE[V_{\alpha, n}(i+\arrival_1,j+\arrival_2)- V_{\alpha, n}(i+\arrival_1, [j-1]^++\arrival_2)] \\
&\leq& \expE[V_{\alpha, n}([i-1]^++\arrival_1,j+\arrival_2)- \\
&&\hspace{0.4cm} V_{\alpha, n}([i-1]^++\arrival_1, [j-1]^++\arrival_2)].
\eean
When $i,j \neq 0$, it is true according to submodularity of $V_{\alpha,n}(i,j)$. Otherwise, both sides of the inequality are 0.

\textbf{Case (ii):} if  $(a^*_{i,j}, a^*_{i+1,j}, a^*_{i,j+1}, a^*_{i+1,j+1})=(0,0,0,0)$, we claim that
\begin{align*}
&\expE[V_{\alpha, n}(i+1+\arrival_1,j+1+\arrival_2)- V_{\alpha, n}(i+1+\arrival_1, j+\arrival_2)]  \\
\leq&\expE[V_{\alpha, n}(i+\arrival_1,j+1+\arrival_2)- V_{\alpha, n}(i+\arrival_1, j+\arrival_2)].
\end{align*}
\annote[AS]{This is obvious from the  submodularity of $V_{\alpha,n}(i,j)$.}{Not clear -what exactly is obvious?}

\textbf{Case (iii):} if $(a^*_{i,j}, a^*_{i+1,j}, a^*_{i,j+1}, a^*_{i+1,j+1})=(0,0,0,1)$, we claim that
\begin{align*}
 & C_T-C_H+ \alpha \expE[V_{\alpha, n}(i+\arrival_1,j+\arrival_2)-\\
 & \hspace{2.4cm} V_{\alpha, n}(i+1+\arrival_1, j+\arrival_2)] \\
\leq& C_H+ \alpha \expE[V_{\alpha, n}(i+\arrival_1,j+1+\arrival_2)-\\
&\hspace{1.5cm} V_{\alpha, n}(i+\arrival_1, j+\arrival_2)].    
\end{align*}
From the submodularity of  $V_{\alpha,n}(i,j)$, it is obtained that
\begin{align*}
&V_{\alpha, n}(i,j)- V_{\alpha, n}(i+1,j)+V_{\alpha, n}(i,j)-V_{\alpha, n}(i,j+1) \\
\leq& V_{\alpha, n}(i,j)- V_{\alpha, n}(i+1,j)+V_{\alpha, n}(i+1,j)-\\
& V_{\alpha, n}(i+1,j+1) \\
=&V_{\alpha, n}(i,j) -V_{\alpha, n}(i+1,j+1).
\end{align*}
Since $a^*_{i+1,j+1}=1$, we have $\Delta \calV_{\alpha, n}(i+1,j+1) \leq 0$, i.e.,
\bean        
&&C_T- 2C_H + \alpha \expE[V_{\alpha, n}(i+\arrival_1,j+\arrival_2)-\\
&&\hspace{2.6cm}V_{\alpha, n}(i+1+\arrival_1,j+1+\arrival_2)] \leq 0.
\eean
\change[AS]{The claim is true because}{The claim follows from the following equation:}
\begin{align*}
&C_T-2C_H+ \alpha \expE[V_{\alpha, n}(i+\arrival_1,j+\arrival_2)- \\
&\hspace{2.5cm}V_{\alpha, n}(i+1+\arrival_1, j+\arrival_2) + \nonumber\\
&\hspace{2.5cm}V_n(i+\arrival_1, j+\arrival_2) - \\
& \hspace{2.5cm}V_{\alpha, n}(i+\arrival_1,j+1+\arrival_2)]  \\
\leq&\Delta \calV_{\alpha,n}(i+1,j+1) \leq 0.
\end{align*}

\textbf{Case (iv):} if  $(a^*_{i,j}, a^*_{i+1,j}, a^*_{i,j+1}, a^*_{i+1,j+1})=(0,0,1,1)$, we claim that
\bean
&& -C_H+ \alpha \expE[V_{\alpha,n}(i+\arrival_1,j+\arrival_2)-\\
&&\hspace{1.7cm}V_{\alpha,n}(i+1+\arrival_1, j+\arrival_2)]  \\
&\leq& C_H([i-1]^+ -i)+ \alpha \expE[V_{\alpha,n}([i-1]^++\arrival_1,j+\arrival_2)-\\
&&\hspace{3.5cm} V_{\alpha,n}(i+\arrival_1, j+\arrival_2)]
\eean
When $i \neq 0$, it is satisfied because $V_{\alpha,n}(i,j)$ is convex. Otherwise, it is true since  $V_{\alpha,n}(i,j)$ is non-decreasing.  

\textbf{Case (v):} if $(a^*_{i,j}, a^*_{i+1,j}, a^*_{i,j+1}, a^*_{i+1,j+1})=(0,1,0,1)$, we claim that
\begin{align*}
& C_H(j - [j-1]^+)+ \alpha \expE[V_{\alpha,n}(i+\arrival_1,j+\arrival_2)-\\
& \hspace{3.5cm} V_{\alpha,n}(i+\arrival_1, [j-1]^++\arrival_2)] \\
\leq& C_H + \alpha \expE[V_{\alpha,n}(i+\arrival_1,j+1+\arrival_2)- V_{\alpha,n}(i+\arrival_1, j+\arrival_2)].
\end{align*}
When $j \neq 0$, it holds since $V_{\alpha,n}(i,j)$ is convex. It is true for other cases because of the non-decreasing $V_{\alpha,n}(i,j)$.

\textbf{Case (vi):} if $(a^*_{i,j}, a^*_{i+1,j}, a^*_{i,j+1}, a^*_{i+1,j+1})=(0,1,1,1)$, we claim that
\begin{align*}
& C_H(j - [j-1]^+)+ \alpha \expE[V_{\alpha,n}(i+\arrival_1,j+\arrival_2)- \\
&\hspace{3.6cm}V_{\alpha,n}(i+\arrival_1, [j-1]^++\arrival_2)] \\
\leq& C_T+ C_H([i-1]^+ - i) + \alpha \expE[V_{\alpha,n}([i-1]^++\arrival_1,j+\arrival_2)-\\
&\hspace{4.4cm} V_{\alpha,n}(i+\arrival_1, j+\arrival_2)].
\end{align*}
Based on the submodularity of $V_{\alpha,n}(i,j)$,  we have
\begin{align*}
&V_{\alpha,n}([i-1]^+,j) -V_{\alpha,n}(i,j) +V_{\alpha,n}(i,[j-1]^+) - V_{\alpha,n}(i,j)\\
\geq& V_{\alpha,n}([i-1]^+,[j-1]^+) -V_{\alpha,n}(i,[j-1]^+) +\\
&V_{\alpha,n}(i,[j-1]^+) - V_{\alpha,n}(i,j)\\
=& V_{\alpha,n}([i-1]^+,[j-1]^+) - V_{\alpha,n}(i,j).
\end{align*}
It is noted that $a^*_{i,j}=0$ and hence  $\Delta \calV_{\alpha,n}(i,j) \geq 0$, i.e.,
\begin{align*}
&C_T + C_H([i-1]^+ +[j-1]^+ -i-j) +\nonumber\\
&\alpha \expE[V_{\alpha,n}([i-1]^+ + \arrival_1, [j-1]^+ + \arrival_1)- \nonumber\\
&V_{\alpha,n}(i+ \arrival_1, j + \arrival_1)] \geq 0.
\end{align*}
Therefore, it can be \change{acquired}{concluded} that 
\begin{align*}
&C_T + C_H([i-1]^+ +[j-1]^+ -i-j) + \\
&\alpha \expE[V_{\alpha,n}([i-1]^++\arrival_1,j+\arrival_2)- \nonumber\\
& \hspace{.5cm}V_{\alpha,n}(i+\arrival_1, j+\arrival_2) +V_{\alpha,n}(i+\arrival_1, [j-1]^++\arrival_2) -\\
& \hspace{.5cm}V_{\alpha,n}(i+\arrival_1,j+\arrival_2)] \\
\geq&  \Delta \calV_{\alpha,n}(i,j) \geq 0.
\end{align*}
\end{proof}

\begin{lemma} \label{lemma:subconvex}
Given $0 < \alpha <1$ and $n \in \mathbb{N}  \cup \{0\}$. If $V_{\alpha,n}(i,j)$ is non-decreasing, submodular, and subconvex, then $V_{\alpha, n+1}(i,j)$ is subconvex.
\end{lemma}

\begin{proof}
We want to show that $V_{\alpha, n+1}(i+1,j+1)- V_{\alpha,  n+1}(i,j) \leq V_{\alpha, n+1}(i+2,j+1)-V_{\alpha,n+1}(i+1,j) $ for all $i$ and $j$. There will be  5 cases of $(a^*_{i,j}, a^*_{i+1,j}, a^*_{i+1,j+1}, a^*_{i+2,j+1})$ that need to be considered. 

\textbf{Case (i):} if  $(a^*_{i,j}, a^*_{i+1,j}, a^*_{i+1,j+1}, a^*_{i+2,j+1})=(1,1,1,1)$, we claim that
\begin{align*}
& C_H(i- [i-1]^+)+ \alpha \expE[V_{\alpha,n}(i+\arrival_1,j+\arrival_2)- \\
&\hspace{3.6cm}V_{\alpha,n}([i-1]^++\arrival_1, [j-1]^++\arrival_2)] \\
\leq& C_H+ \alpha \expE[V_{\alpha,n}(i+1+\arrival_1,j+\arrival_2)- \\
&\hspace{1.6cm} V_{\alpha,n}(i+\arrival_1, [j-1]^++\arrival_2)].
\end{align*}
When $i,j \neq 0$, it is true according to the subconvexity of $V_{\alpha,n}(i,j)$. The argument is satisfied for $i=0, j \neq 0$ due to the the non-decreasing $V_{\alpha,n}(i,j)$, and for the case $i \neq 0, j=0$ due to the convexity of $V_{\alpha,n}(i,j)$. Otherwise, it holds according to the non-decreasing property.

\textbf{Case (ii):} if $(a^*_{i,j}, a^*_{i+1,j}, a^*_{i+1,j+1}, a^*_{i+2,j+1})=(0,0,0,0)$, we claim that
\begin{align*}
&\expE[V_{\alpha,n}(i+1+\arrival_1,j+1+\arrival_2)- V_{\alpha,n}(i+\arrival_1, j+\arrival_2)] \\
\leq& \expE[V_{\alpha,n}(i+2+\arrival_1,j+1+\arrival_2)-\\
&\hspace{.3cm}V_{\alpha,n}(i+1+\arrival_1, j+\arrival_2)]. 
\end{align*}
 The above results from the subconvexity of $V_{\alpha,n}(i,j)$.
 
\textbf{Case (iii):} if $(a^*_{i,j}, a^*_{i+1,j}, a^*_{i+1,j+1}, a^*_{i+2,j+1})=(0,0,0,1)$, we claim that
\begin{align*}
&2C_H+ \alpha \expE[V_{\alpha,n}(i+1+\arrival_1,j+1+\arrival_2)- \\
&\hspace{1.6cm}V_{\alpha,n}(i+\arrival_1, j+\arrival_2)] \leq C_T. 
\end{align*}
Since $a^*_{i+1,j+1}=0$, we have $\Delta \calV_{\alpha,n}(i+1,j+1) \geq 0$, i.e.,
\begin{align*}
    & C_T- 2C_H + \alpha \expE[V_{\alpha,n}(i+\arrival_1,j+\arrival_2)-\\
     &\hspace{2.5cm}V_{\alpha,n}(i+1+\arrival_1,j+1+\arrival_2)] \geq 0.
\end{align*}
Hence the claim is verified.

\textbf{Case (iv):} if  $(a^*_{i,j}, a^*_{i+1,j}, a^*_{i+1,j+1}, a^*_{i+2,j+1})=(0,0,1,1)$, it is trivial since the both  $V_{\alpha, n+1}(i+1,j+1)- V_{\alpha,  n+1}(i,j)$ and $V_{\alpha, n+1}(i+2,j+1)-V_{\alpha,n+1}(i+1,j) $  are zeros.

\textbf{Case (v):} if $(a^*_{i,j}, a^*_{i+1,j}, a^*_{i+1,j+1}, a^*_{i+2,j+1})=(0,1,1,1)$, we claim that
\bean
C_T &\leq& C_H(1+ j- [j-1]^+)+ \\
&&\alpha \expE[V_{\alpha,n}(i+1+\arrival_1,j+\arrival_2)-\\
&&\hspace{.5cm} V_{\alpha,n}(i+\arrival_1, [j-1]^++\arrival_2)]. 
\eean
Notice that $a^*_{i+1,j}=1$, so $\Delta \calV_{\alpha,n}(i+1,j) \leq 0$, i.e.,
\bean
        &&C_T- C_H(1+ j- [j-1]^+) +\nonumber\\
        &&\alpha \expE[V_{\alpha,n}(i+\arrival_1,[j-1]^++\arrival_2)-\\
        &&\hspace{.5cm} V_{\alpha,n}(i+1+\arrival_1,j+\arrival_2)] \leq 0.
 \eean
\end{proof}

Based on the properties of $V_{\alpha}(i,j)$, we are ready to state the optimality of the threshold type policy in terms of the total expected discounted cost.   
\begin{theorem}
\label{thm:thresold-discounted} 
For the MDP$\{(Q_t,a_t), t \geq 0\}$ with any i.i.d. arrival processes to both queues, there exists an $\alpha$-optimal policy that is of threshold type. Given $Q^{(2)}_t$, the $\alpha$-optimal policy is monotone w.r.t. $Q^{(1)}_t$, and vice versa.
\end{theorem}
\begin{proof}
We prove by induction. $V_{\alpha,0}(i,j)=0$ is non-decreasing, submodular, and subconvex, that leads to the non-decreasing $\min \{a' \in \argmin_{a \in \{0,1\}} \calV_{\alpha,0}(i,j,a)\}$ based on Corollary \ref{cor:itr-monotone-policy}. These properties propagate as $n$ goes to infinity according to lemmas \ref{lemma:increasing function}, \ref{lemma:submodular}, \ref{lemma:subconvex}, and Corollary~\ref{cor:itr-monotone-policy}. 
\end{proof}

Thus far, the $\alpha$-optimal policy is characterized.  A useful relation between the average-optimal policy and the $\alpha$-optimal policy is described in the following lemma. 

\begin{lemma}[\cite{stationary-policy-Sennott}, Lemma and Theorem (i)]
\label{lemma:deterministic-policy-cond}
Consider  MDP$\{(Q_t, a_t), t \geq 0\}$. Let $\{ \alpha_n\}$ converging to 1 be any sequence of discount factors associated with the $\alpha$-optimal policy $\{\theta_{\alpha_n}(i,j)\}$. There exists a subsequence $\{\beta_n\}$  and a stationary policy $\theta^*(i,j)$ that is the limit point of  $\{\theta_{\beta_n}(i,j)\}$. If the three conditions in Lemma \ref{lemma:stationary-policy-exist} are satisfied, $\theta^*(i,j)$ is the average-optimal policy for Eq. (\ref{eq:long-time-average-cost}).
\end{lemma}

\begin{theorem} \label{thm:thresold}
Consider any i.i.d. arrival processes to both queues. For the MDP$\{(Q_t, a_t), t \ge 0\}$, the average-optimal policy is of threshold type. There exist the optimal thresholds $L^*_1$ and $L^*_2$ so that the optimal deterministic action in states $(i,0)$ is to wait if $i \leq L^*_1$, and to transmit without coding if $i > L^*_1$; while in state $(0,j)$ is to wait if $j \leq L^*_2$, and to transmit without coding if $j > L^*_2$. 
\end{theorem}
\begin{proof}
Let $(\tilde{i},0)$ be any state which average-optimal policy is to transmit, i.e., $\theta^*(\tilde{i},0)=1$ in Lemma \ref{lemma:deterministic-policy-cond}. Since there is a sequence of discount factors $\{\beta_{n}\}$ such that $\theta_{\beta_n}(i,j) \rightarrow \theta^*(i,j)$, then there exists $N>0$ so that $\theta_{\beta_n}(\tilde{i},0)=1$ for all $n \geq N$. Due to the monotonicity of $\alpha$-optimal policy in Theorem~\ref{thm:thresold-discounted}, $\theta_{\beta_n}(i,0)=1$ for all $i \geq \tilde{i}$ and $n \geq N$. Therefore, $\theta^*(i,0)=1$ for all $i \geq \tilde{i}$. To conclude, the average-optimal policy is of threshold type.
\end{proof}

\section{Obtaining the Optimal Deterministic Stationary Policy}

We have shown in the previous sections that the average-optimal policy is stationary, deterministic and of threshold type, so we only need to consider the subset of deterministic stationary policies. Given the thresholds of the both queues, the MDP is reduced to a Markov chain. The next step is to find the optimal threshold. First note that the condition $\expE[\arrival_i] < 1$  might  not be  sufficient for the stability of the queues since the threshold based policy leads to an average service rate lower than 1 packet per time slot. In the following theorem, we claim that the conditions $\expE[\arrival_i^2] < \infty$ and $\expE[\arrival_i] < 1$ for  $i=1,2$ are enough for the stability of the queues. 

\begin{theorem}
For the MDP$\{(Q_t, a_t), t \ge 0\}$ with $\expE[\arrival_i^2] < \infty$ and $\expE[\arrival_i] < 1$ for  $i=1,2$. The reduced Markov chain from applying the stationary and deterministic threshold based policy to MDP is positive recurrent, i.e., the stationary distribution exists. 
\end{theorem}

\begin{proof}
The proof is based on Foster-Lyapunov theorem \cite{foster-lyapunov} associated with the Lyapunov function $\mathcal{L}(x,y)=x^2+y^2$. Notice that
\bean
Q^{(i)}_{t+1}&=&[Q^{(i)}_t-a_t]^++ \arrival_i \\
&=& Q^{(i)}_t-a_t +U^{(i)}_t+\arrival_i,
\eean
where
\bean
U^{(i)}_t= \left\{
\begin{array}{ll}
0 & \text{if\,\,} Q^{(i)}_t -a_t \geq 0  \\ 
1 & \text{if\,\,} Q^{(i)}_t -a_t =-1. \\ 
\end{array}
\right.
\eean
Then it can be observed that  
\begin{align}
&\expE \left[ \mathcal{L}(Q^{(1)}_{t+1}, Q^{(2)}_{t+1}) - \mathcal{L}(Q^{(1)}_{t}, Q^{(2)}_{t}) | Q^{(1)}_{t}=x, Q^{(2)}_{t}=y \right] \nonumber\\
=& \expE \left[ \sum^{2}_{i=1} (Q^{(i)}_t -a_t+ U^{(i)}_t + \arrival_i)^2  | Q^{(1)}_{t}=x, Q^{(2)}_{t}=y  \right] -\nonumber\\
&(x^2+y^2) \nonumber
\end{align}
\begin{align}
=& \sum^{2}_{i=1}  \expE \left[ (Q^{(i)}_t -a_t+  \arrival_i)^2  | Q^{(1)}_{t}=x, Q^{(2)}_{t}=y  \right] +  \nonumber\\
&\sum^{2}_{i=1} \expE \left[ (U^{(i)}_t )^2  | Q^{(1)}_{t}=x, Q^{(2)}_{t}=y  \right]+ \nonumber\\ 
&\sum^{2}_{i=1} \expE \left[ 2 U^{(i)}_t (Q^{(i)}_t -a_t+  \arrival_i)  | Q^{(1)}_{t}=x, Q^{(2)}_{t}=y  \right] - \nonumber\\
& (x^2+y^2) \nonumber\\
\leq& \sum^{2}_{i=1}  \expE \left[ (Q^{(i)}_t -a_t+  \arrival_i)^2  | Q^{(1)}_{t}=x, Q^{(2)}_{t}=y  \right] +2 +\nonumber\\
& 2\expE[\arrival_1] + 2\expE[\arrival_2] - (x^2+y^2) \label{eq:stable_inequality_1} \\
=& 2 x \expE \left[ \arrival_1- a_t | Q^{(1)}_{t}=x, Q^{(2)}_{t}=y \right] + \nonumber\\
&2 y \expE \left[ \arrival_2- a_t | Q^{(1)}_{t}=x, Q^{(2)}_{t}=y \right] +\nonumber\\
&  \expE \left[ (\arrival_1- a_t)^2 | Q^{(1)}_{t}=x, Q^{(2)}_{t}=y \right] +\nonumber\\
&  \expE \left[( \arrival_2- a_t)^2 | Q^{(1)}_{t}=x, Q^{(2)}_{t}=y \right] +\nonumber\\
& 2\expE[\arrival_1] +  2\expE[\arrival_2] +2 \nonumber\\
\leq& 2 x \expE \left[ \arrival_1- a_t | Q^{(1)}_{t}=x, Q^{(2)}_{t}=y \right] +\nonumber\\
& 2 y \expE \left[ \arrival_2- a_t | Q^{(1)}_{t}=x, Q^{(2)}_{t}=y \right] +\nonumber\\
& \expE[\arrival^2_1] +1 + \expE[\arrival^2_2]+1 + 2\expE[\arrival_1] +  2\expE[\arrival_2] +2 \label{eq:stable_inequality_2}\\
=& \expE[\arrival^2_1]+\expE[\arrival^2_2]+2\expE[\arrival_1] +  2\expE[\arrival_2] +4+\nonumber\\
&\left\{
\begin{array}{ll}
2x(\expE[\arrival_1]-1)+2y(\expE[\arrival_2]-1) & \text{if}\,\, (x,y) \in \mathcal{B}^c\\
2x\expE[\arrival_1] +2y\expE[\arrival_2]  & \text{if\,\,} (x,y) \in \mathcal{B},
\end{array} 
\right. \label{eq:stable_inequlity_3}
\end{align}
where $\mathcal{B}=\{(x,y): (x=0, y \leq L_2) \text{\,\,or\,\,} (x \leq L_1, y=0) \}$. The inequality (\ref{eq:stable_inequality_1}) comes from $(U^{(i)}_t)^2 \leq 1$ and $(Q^{(i)}_t-a_t) U^{(i)}_t \leq 0$, while $\expE[a_t] \leq 1$ results in Eq. (\ref{eq:stable_inequality_2}). Since  $\expE[\arrival_i^2] < \infty$ and $\expE[\arrival_i] < 1$ for  $i=1, 2$,  the value in Eq.  (\ref{eq:stable_inequlity_3}) is negative for $(x,y) \in \mathcal{B}^c$ and is bounded for $(x,y) \in \mathcal{B}$. Then the result immediately follows from Foster-Lyapunov theorem. 
\end{proof}

We realize that if $\expE[\arrival_i^2] < \infty$ and $\expE[\arrival_i] < 1$ for  $i =1, 2$, then there exists a stationary threshold type policy that is average-optimal and can be obtained from the reduced Markov chain. The following theorem gives an example of how to compute the optimal thresholds. 

\begin{theorem}
Consider the Bernoulli arrival process. The optimal thresholds $L^*_1$ and $L^*_2$ are
\bean
(L^*_1,L^*_2) = \argmin_{L_1,L_2}C_T \mathcal{T}(L_1,L_2) + C_H \mathcal{H}(L_1,L_2),
\eean 
where
\bean
\mathcal{T}(L_1,L_2) &=& p^{(1)}_1 p^{(2)}_1 \pi_{0,0} + p^{(2)}_1 \sum_{i=1}^{L_1} \pi_{i,0} + p^{(1)}_1 \sum_{j=1}^{L_2} \pi_{0,j} +\\
&&p^{(1)}_1 p^{(2)}_0 \pi_{L_1,0} +p^{(1)}_0 p^{(2)}_1 \pi_{0,L_2}; \\
\mathcal{H}(L_1,L_2) &=& \sum_{i=1}^{L_1} i \pi_{i,0} + \sum_{j=1}^{L_2} j \pi_{0,j},
\eean
for which
\bean
\pi_{0,0} &=& \frac{1}{\left( \frac{1 - \zeta^{L_1+1}}{1 - \zeta} \right) + \left( \frac{1 - 1/\zeta^{L_2+1}}{1 - 1/\zeta} \right) - 1}; \\
\pi_{i,0} &=&\zeta^i  \pi_{0,0}; \\
\pi_{0,j} &=& \pi_{0,0}/\zeta^j;  \\
\zeta &=& \frac{p^{(1)}_1 p^{(2)}_0}{p^{(1)}_0 p^{(2)}_1}. 
\eean
\label{thm:optimal-threshold}
\end{theorem}

\begin{proof} 
Let $Y^{(i)}_t$ be the number of type $i$ packets at the $t^{\text{th}}$ slot \textit{after} transmission. It is crucial to note that this observation time is different from when the MDP is observed. Then the bivariate stochastic process $\{(Y^{(1)}_t,Y^{(2)}_t), t \ge 0\}$ is a discrete-time Markov chain which state space is smaller than the original MDP, i.e.,  $(0,0)$, $(1,0)$, $(2,0)$, $\cdots$, $(L_1,0)$, $(0,1)$, $(0,2)$, $\cdots$, $(0,L_2)$.
Define $\zeta$ as a parameter such that
\bean
\zeta = \frac{p^{(1)}_1 p^{(2)}_0}{p^{(1)}_0 p^{(2)}_1}. 
\eean
Then, the balance equations for $0 < i \leq L_1$ and $0 < j \leq L_2$ are:
\bean
\pi_{i,0} &=& \zeta  \pi_{i-1,0};  \\
\zeta  \pi_{0,j} &=& \pi_{0,j-1}.
\eean
Since $\pi_{0,0} + \sum_{i,j} \pi_{i,0} + \pi_{0,j} = 1$, we have
\bean
\pi_{0,0} &=& \frac{1}{\left( \frac{1 - \zeta^{L_1+1}}{1 - \zeta} \right) + \left( \frac{1 - 1/\zeta^{L_2+1}}{1 - 1/\zeta} \right) - 1}.
\eean

The expected number of transmissions per slot is
\bean
\mathcal{T}(L_1,L_2) &=& p^{(1)}_1 p^{(2)}_1 \pi_{0,0} + p^{(2)}_1 \sum_{i=1}^{L_1} \pi_{i,0} + p^{(1)}_1 \sum_{j=1}^{L_2} \pi_{0,j} +\\
&&p^{(1)}_1 p^{(2)}_0 \pi_{L_1,0} +  p^{(1)}_0 p^{(2)}_1 \pi_{0,L_2}.
\eean
The average number of packets in the system at the beginning of each slot is
\bean
\mathcal{H}(L_1,L_2) = \sum_{i=1}^{L_1} i \pi_{i,0} + \sum_{j=1}^{L_2} j \pi_{0,j}.
\eean
Thus upon minimizing we get the optimal thresholds $L^*_1$ and $L^*_2$.
\end{proof}

Whenever $C_H > 0$, it is relatively straightforward to obtain $L^*_1$ and $L^*_2$. Since it costs $C_T$ to transmit a packet and $C_H$ for a packet to wait for a slot, it would be better to transmit a packet than make a packet wait for more than $C_T/C_H$ slots. Thus $L^*_1$ and $L^*_2$ would always be less than $C_T/C_H$. Hence, by completely
enumerating between $0$ and $C_T/C_H$ for both $L_1$ and $L_2$, we can obtain $L^*_1$ and $L^*_2$. One could perhaps find faster techniques than complete enumeration, but it certainly serves the purpose.

Subsequently, we study a special case, $p^{(1)}_1=p^{(2)}_1 \triangleq p$, in Theorem \ref{thm:optimal-threshold}. Then $L_1=L_2 \triangleq L$ as  both arrival processes are identical. It can be calculated that $\zeta=1$ and $\pi_{i,j}=1/(2L+1)$ for all $i,j$, and 
\bean
\mathcal{T}(L)= \frac{2pL+2p-p^2}{2L+1};
\eean
\bean
\mathcal{H}(L)= \frac{L^2+L}{2L+1}.
\eean
Define $\nu=C_T/C_H$. The optimal threshold is 
\bean
L^*(p,\nu)= \argmin_L \frac{\nu(2pL+2p-p^2)+L+L^2}{2L+1}.
\eean
By taking the derivative, we obtain that $L^*=0$ if $\nu < 1/(2p-2p^2)$ and otherwise, 
\bean
L^*(p,\nu)=\frac{-1+\sqrt{1-2(1-2 \nu p+2 \nu p^2)}}{2}.
\eean

We can observe that $L^*(p,\nu)$ is a concave function w.r.t. $p$. Given $\nu$ fixed, $L^*(1/2,\nu) = (\sqrt{\nu-1}-1)/2$ is the largest optimal threshold among various values of $p$. When $p < 1/2$, the optimal-threshold decreases as there is a relatively lower probability for packets in one queue to wait for a coding pair in another queue. When $p >1/2$, there will be a coding pair already in the relay node with a higher probability, and therefore the optimal-threshold also decreases. Moreover, $L^*(1/2,\nu)=\mathcal{O}(\sqrt{\nu})$, so the maximum optimal threshold grows with the square root of $\nu$, but not linearly. When $p$ is very small, $L^*(p,\nu)=\mathcal{O}(\sqrt{\nu p})$ grows slower than $L^*(1/2,\nu)$.

%% file: numericals.tex

\section{Numerical Studies}\label{sec:numerical-1}

In this section we present several numerical results to compare the performance of different policies in the single relay setting as well as in the line network. 
We analyzed the following policies:

\begin{enumerate}
\item Opportunistic Coding (OC):  this policy does not waste any opportunities for transmitting the packets. That is when a packet arrives, coding is performed if a coding opportunity exists, otherwise transmission takes place immediately.
\item Queue-length  based threshold (QLT):  this stationary deterministic policy applies the thresholds, proposed by Theorem~\ref{thm:optimal-threshold}, on the queue lengths. 
\item Queue-length-plus-Waiting-time-based (QL+WT) thre-sholds: this is a history dependent policy which takes into account the waiting time of the packets in the queues as well as the queue lengths. That is a packet will be transmitted (without coding), if the queue length hits the threshold or the head-of-queue packet has been waiting for at least some pre-determined amount of time. The optimal waiting-time thresholds are found using exhaustive search through stochastic simulations for the given arrival distributions.
\item Waiting-time (WT) based threshold: this is another history dependent policy that \emph{only} considers the waiting times of the packets,  in order to schedule the transmissions. The optimum waiting times of the packets are found through exhaustive search. 
\end{enumerate}

\begin{figure}[h]
\begin{center}
\includegraphics[width=0.5\textwidth]{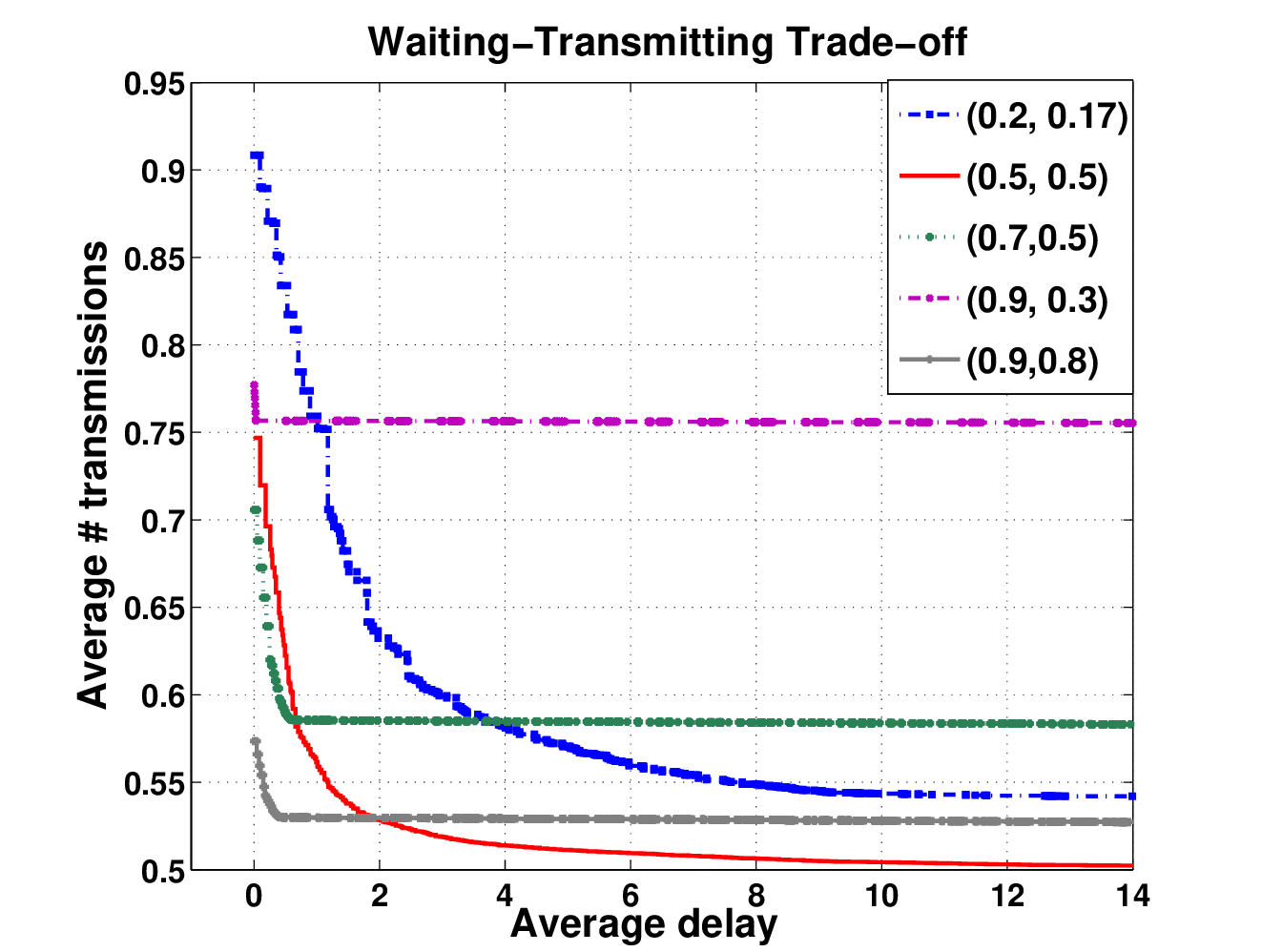}
\caption{Trade-off between average delay and number of transmissions in a single relay using queue-length based threshold (QLT) policy for different Bernoulli arrival rates $(p_1,p_2)$.}
\label{Fig:tradeoff}
\end{center}
\end{figure}

\begin{figure}[h]
\begin{center}
\includegraphics[width=0.5\textwidth]{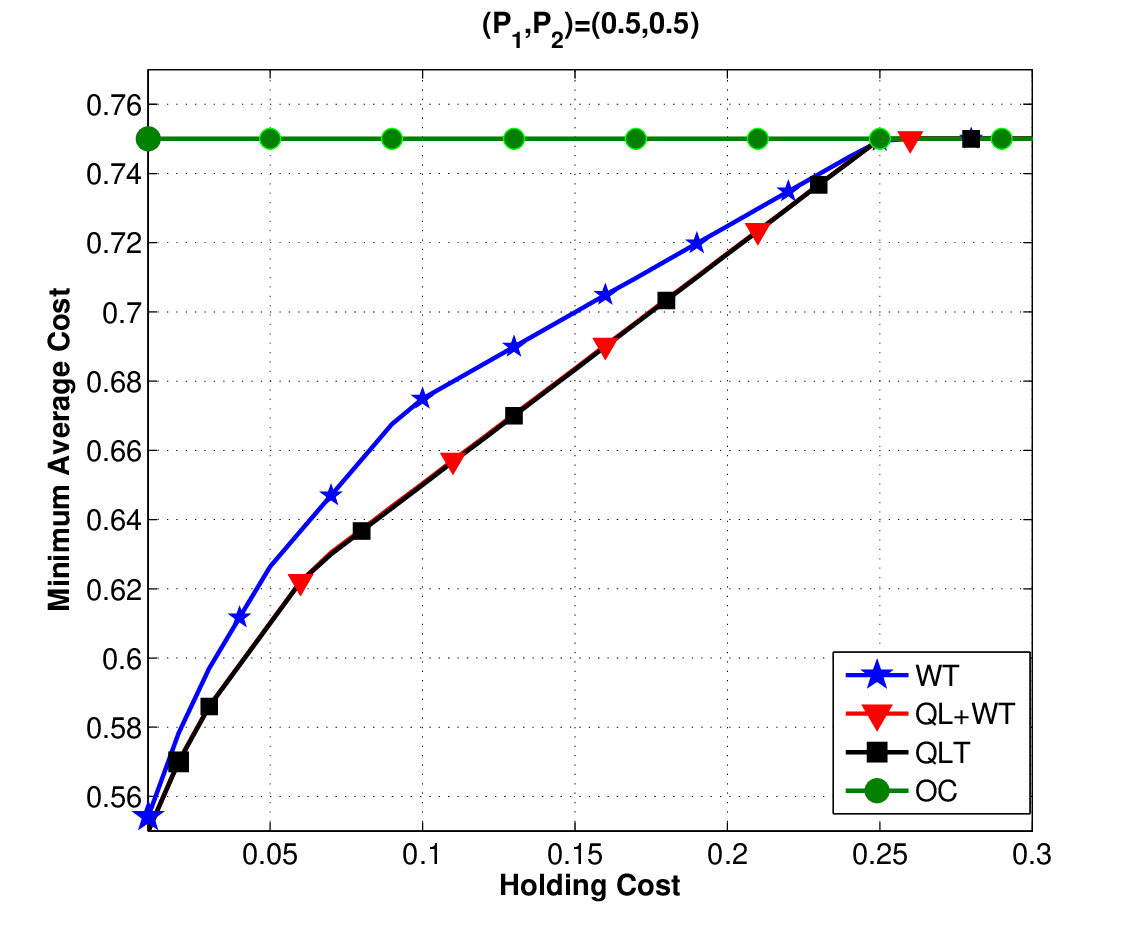}
\caption{Comparison of the minimum average cost (per packet) in a single relay with Bernoulli arrival rates $(0.5,0.5)$, for different policies, where the costs are normalized by the transmission cost.}
\label{Fig:BER55}
\end{center}
\end{figure}

\begin{figure}[h]
\begin{center}
\includegraphics[width=0.5\textwidth]{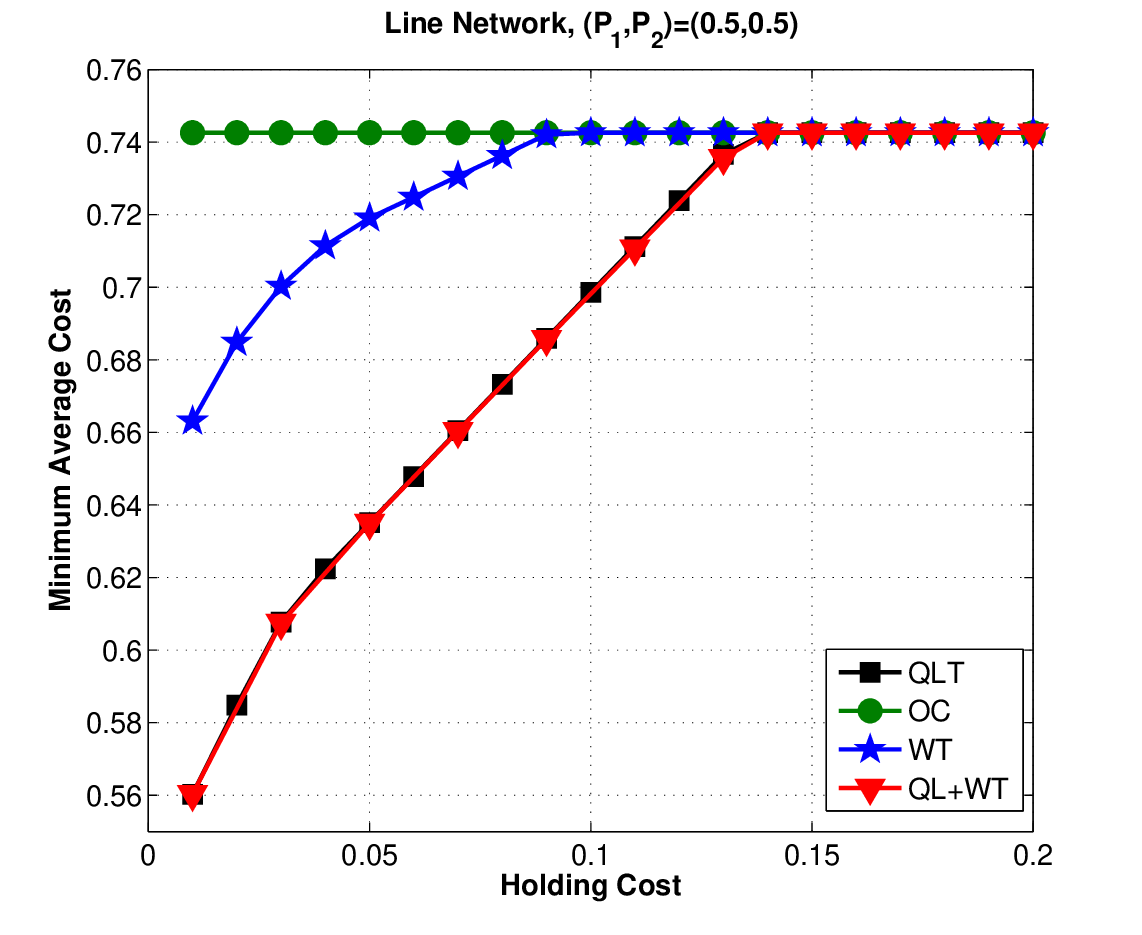}
\caption{Comparison of different policies in a line network with two intermediate nodes and two Bernoulli flows with mean arrival rates $(0.5,0.5)$.}
\label{Fig:LINE55}
\end{center}
\end{figure}

We simulate these policies on two different cases: (i) the single relay network with Bernoulli arrivals (Figures \ref{Fig:tradeoff} and \ref{Fig:BER55}) and (ii) a line network with $4$ nodes, in which the sources are Bernoulli (Figure \ref{Fig:LINE55}). Note that in case (ii), since the departures from one queue determine the arrivals into the other queue, the arrival processes are significantly different from Bernoulli.
Our simulations are done in Java and for each scenario we report the average results of $10^5$ iterations.

As expected, for the single relay network, the QLT policy has the optimal performance and the QL+WT policy does not have any advantage. Our simulation results indicate that QLT policy also exhibits a near optimal performance for the line network. We also observe, from the simulation results for the waiting-time-based policy, that making decisions based on waiting time alone leads to a suboptimal performance. In all experiments, the opportunistic policy has the worst possible performance. 

The results are intriguing as they suggest that achieving a near-perfect trade-off between waiting and transmission costs is possible using simple policies; moreover, coupled with optimal network-coding aware routing policies like the one in our earlier work \cite{population-game-shakkottai}, have the potential to exploit the positive externalities that network coding offers.

%% file: extension.tex

\section{Extensions}\label{sec:extension}

We have seen that the average-optimal policy is stationary and threshold based for the i.i.d. arrival process with the service rate of 1 packet per time slot. Two more general models are discussed here. We focus on the character of the optimality equation which results in the structure of the average-optimal policy. 

\subsection{Batched service}
Assume that  the relay $R$ can serve a group of packets with the size of $\mathcal{M}$ at end of the time slot.  At the end of every time slot, relay $R$ decides to transmit, $a_t=1$, or to wait $a_t=0$. The holding cost per unit time for a packet is $C_H$, while $C_T$ is the cost to transmit a  batched packet. Then the immediate cost is 
\bean
&&C^{(\mathcal{M})}(Q_t, a_t) \\
&=& C_H ([Q^{(1)}_t-a_t \mathcal{M}]^+ + [Q^{(2)}_t-a_t  \mathcal{M}]^+) + C_T a_t.
\label{eq:immediate cost-k}
\eean
We also want to find the optimal policy $\theta^*$ that minimizes the long-time average cost $V^{(\mathcal{M})}(\theta)$, called $\mathcal{M}$-MDP$\{(Q_t, a_t), t \geq 0\}$ problem,
\bean
V^{(\mathcal{M})}(\theta) = \lim_{K \rightarrow \infty} \frac{1}{K+1} \expE_{\theta} \left[ \sum_{t=0}^K C^{(\mathcal{M})}(Q_t,a_t) | Q_0 = (0,0) \right]. \label{eq:long-time-average-cost-k}
\eean
Notice that  the best policy might not just transmit when both queues are non-empty. When $\mathcal{M} > 1$, $R$ might also want to wait even if $Q^{(1)}_t Q^{(2)}_t >0$ because the batched service of size less than $\mathcal{M}$ has the same transmission cost $C_T$. The optimality equation of the expected $\alpha$-discounted cost is revised as
\begin{align*}
V^{(\mathcal{M})}_{\alpha}(i,j) = &\min_{a \in \{0,1\}} \Bigl[ C_H ([i-a \mathcal{M}]^+ + [j-a \mathcal{M}]^+) + C_T a + \nonumber\\
& \hspace{.2cm}\expE[V^{(\mathcal{M})}_{\alpha}([i-a \mathcal{M}]^++\arrival_1,[j-a \mathcal{M}]^++ \arrival_2)] \Bigr].
 \label{eq:optimality-eq-k}
\end{align*}
We can get the following results.
\begin{theorem}
Given $\alpha$ and $\mathcal{M}$, $V^{(\mathcal{M})}_{\alpha}(i,j)$ is non-decreasing, submodular, and $\mathcal{M}$-subconvex. Moreover, there is an $\alpha$-optimal policy that is of threshold type.  Fixed $j$, the $\alpha$-optimal policy is monotone w.r.t. $i$, and vice versa.
\end{theorem}
\begin{theorem} \label{thm:thresold-new}
Consider any i.i.d. arrival processes to both queues. For the $\mathcal{M}$-MDP$\{(Q_t, a_t), t \ge 0\}$, the average-optimal policy is of threshold type.  Given $j=\tilde{j}$ fixed, there exists the optimal threshold $L^*_{\tilde{j}}$ such that the optimal stationary and deterministic policy in state $(i, \tilde{j})$ is to wait if $i \leq L^*_{\tilde{j}}$, and to transmit if $i > L^*_{\tilde{j}}$. Similar argument \change{is true}{holds} for the other queue. 
\end{theorem}

\subsection{Markov-Modulated Arrival Process}
While the i.i.d. arrival process is examined so far, a specific arrival process with memory is studied here, i.e., Markov-modulated arrival process (MMAP). The service capacity of $R$ is focused on $\mathcal{M}=1$ packet.  Let  $\mathcal{N}^{(i)}=\{0, 1, \cdots, N^{(i)}\}$ be the state space of MMAP at node $i$, with the transition probability $p^{(i)}_{k,l}$ where $k, l \in \mathcal{N}^{(i)}$. Then the number of packets generated by the node $i$ at time $t$ is $\mathcal{N}^{(i)}_t \in \mathcal{N}^{(i)}$. Then the  decision of $R$ is made based on the observation of $(Q^{(1)}_t, Q^{(2)}_t, \mathcal{N}^{(1)}_t, \mathcal{N}^{(2)}_t)$. Similarly, the objective is to find the optimal policy that minimizes the  long-term average cost, named MMAP-MDP$\{((Q^{(1)}_t, Q^{(2)}_t, \mathcal{N}^{(1)}_t, \mathcal{N}^{(2)}_t), a_t): t \geq 0\}$ problem. The optimality equation of the expected $\alpha$-discounted cost becomes
\begin{align*}
&V^{\text{MMAP}}_{\alpha}(i,j,n_1,n_2) \\
 =&\min_{a \in \{0,1\}} [ C_H ([i-a]^+ + [j-a]^+) + C_T a + \nonumber\\
& \hspace{.3cm} \alpha \sum_{k=0}^{N^{(1)}} \sum_{l=0}^{N^{(2)}} p^{(1)}_{n_1,k} p^{(2)}_{n_2,l} V^{\text{MMAP}}_{\alpha}([i-a]^++k, [j-a]^++l,k,l)].
\end{align*}
Then we can conclude the following results.
\begin{theorem}
Given $n_1 \in \mathcal{N}^{(1)}$ and $n_2 \in \mathcal{N}^{(2)}$, $V^{\text{MMAP}}_{\alpha}(i,j,n_1,n_2)$ is non-decreasing, submodular, and subconvex w.r.t. $i$ and $j$. Moreover, there is an $\alpha$-optimal policy that is of threshold type.  Fixed $n_1$ and $n_2$, the $\alpha$-optimal policy is monotone w.r.t. $i$ when $j$ is fixed, and vice versa.
\end{theorem}
\begin{theorem} \label{thm:thresold-new2}
Consider any MMAP arrival process. For the MMAP-MDP$\{((Q^{(1)}_t, Q^{(2)}_t, \mathcal{N}^{(1)}_t, \mathcal{N}^{(2)}_t), a_t): t \geq 0\}$, the  average-optimal policy is of multiple thresholds type.   There exists a set of optimal thresholds $\{L^*_{1,n_1,n_2}\}$ and $\{L^*_{2,n_1,n_2}\}$, where $n_1 \in \mathcal{N}^{(1)}$ and $n_2 \in \mathcal{N}^{(2)}$, so that the optimal stationary decision in states $(i,0,n_1,n_2)$ is to wait if $i \leq L^*_{1,n_1,n_2}$, and to transmit without coding if $i >L^*_{1,n_1,n_2}$; while in state $(0,j,n_1,n_2)$ is to wait if $j \leq L^*_{2,n_1,n_2}$, and to transmit without coding if $j > L^*_{2,n_1,n_2}$. 
\end{theorem}

%% file: conclusion.tex
\section{Conclusion}

In this paper we investigate the delicate trade-off between waiting and transmitting using network coding. We started with the idea of exploring the whole space of history dependent policies, but showed step-by-step how we could move to simpler regimes, finally culminating in a stationary deterministic queue-length threshold based policy. The policy is attractive because its simplicity enables us to characterize the thresholds completely, and we can easily illustrate its performance on multiple networks. We showed by simulation how the performance of the policy is optimal in the Bernoulli arrival scenario, and how it also does well in other situations such as for line networks. Our results also have some bearing on the general problem of queuing networks with shared resources that we will explore in the future.

\section{Acknowledgement}
The authors are grateful to Rajesh Sundaresan, Vivek Borkar, and Daren B. H. Cline for the useful discussions. This material is based upon work partially supported by the AFOSR under contract No. FA9550-13-1-0008.